\input ./style/arxiv-general.cfg
\documentclass[aos,MSNbibl,dvips]{arximspdf}
\makeatletter
   \@ifpackageloaded{graphicx}{}{\usepackage{graphicx}}
\makeatother
\usepackage{mathbh}

\doi{10.1214/15-AOS1316}
\volume{43}
\issue{4}
\pubyear{2015}
\firstpage{1568}
\lastpage{1595}
\docsubty{FLA}

\makeatletter
\newcommand{\rrvert}{\vert}
\newcommand{\llvert}{\vert}
\newtheorem{teo}{Theorem}
\newproclaim{defi}{Definition}
\newproclaim{qua}{Quadrant}
\newtheorem{lem}{Lemma}
\newtheorem{prop}{Proposition}
\newcommand{\ee}{e}
\newcommand{\re}{\mathbb{R}}
\newcommand{\simdist}{\stackrel{\mathcal{D}}{\sim}}
\newcommand{\likely}{\mathcal{L}}
\newcommand{\convdist}{\stackrel{\mathcal{D}}{\rightarrow}}
\newcommand{\xo}{x_0}
\newcommand{\za}[1]{\stackrel{a}{#1}}
\newcommand{\zb}[1]{\stackrel{b}{#1}}
\newcommand{\zc}[1]{\stackrel{c}{#1}}
\newcommand{\zd}[1]{\stackrel{d}{#1}}
\makeatother

\begin{document}
\begin{frontmatter}

\title{Robustness to outliers in location--scale parameter model using
log-regularly varying distributions}
\runtitle{Robustness in location--scale parameter model}

\begin{aug}
\author[A]{\fnms{Alain}~\snm{Desgagn\'{e}}\corref{}\ead[label=e1]{desgagne.alain@uqam.ca}}
\runauthor{A. Desgagn\'{e}}
\affiliation{Universit\'{e} du Qu\'{e}bec \`{a} Montr\'{e}al}
\address[A]{D\'{e}partement de math\'{e}matiques\\
Universit\'{e} du Qu\'{e}bec \`{a} Montr\'{e}al\\
C.P. 8888, Succursale Centre-ville\\
Montr\'{e}al, Qu\'{e}bec\\
Canada H3C 3P8\\
\printead{e1}}
\end{aug}

%
\received{\smonth{8} \syear{2014}}
%
\revised{\smonth{1} \syear{2015}}

%
\begin{abstract}
Estimating the location and scale parameters is common in statistics,
using, for instance, the well-known sample mean and standard deviation.
However, inference can be contaminated by the presence of outliers if
modeling is done with light-tailed distributions such as the normal
distribution.
In this paper, we study robustness to outliers in location--scale
parameter models using both the Bayesian and frequentist approaches.
We find sufficient conditions (e.g., on tail behavior of the model) to
obtain whole robustness to outliers,
in the sense that the impact of the outliers gradually decreases to
nothing as the conflict grows infinitely.
To this end, we introduce the family of log-Pareto-tailed symmetric
distributions that belongs to the larger family of log-regularly
varying distributions.
\end{abstract}

%
\end{frontmatter}

\section{Introduction}\label{sec-intro}

In Bayesian analysis, outlying observations and prior misspecification
may contaminate the posterior inference.
For instance, a group of
observations may suggest a quite different posterior inference than
that proposed by the prior and the rest of data.
Using light-tailed distributions such as the normal can lead to an
undesirable compromise where
the posterior distribution concentrates on an area that is not
supported by any sources of information.
The conflict is usually
resolved automatically by modeling with heavy-tailed distributions, in
favor of the sources of information
with the lightest tails.
O'Hagan and Pericchi \cite{r16} refer to this situation as the \textit
{theory of conflict resolution in Bayesian statistics}, in their
extensive review of the literature on that topic.

Conflict resolution in Bayesian analysis was first described by De
Finetti \cite{r7}.
The theory has mostly been developed for location parameter inference;
see, for instance,
Dawid \cite{r6}; O'Hagan \cite{r13,r14,r15}; Angers \cite{r5};
Desgagn\'e and Angers \cite{r10};
Kumar and Magnus \cite{r12}; Andrade and Omey \cite{r4}; Andrade,
Dorea and Guevara Otiniano \cite{r1}.

The theory on pure scale parameter inference was first analyzed by
Andrade and O'Hagan \cite{r2},
who considered partial robustness using regularly varying distributions
(see also Andrade and Omey \cite{r4} and
Andrade, Dorea and Guevara Otiniano \cite{r1}, who generalize their
work of partial robustness), and then
by Desgagn\'e \cite{r8}, who considered whole robustness using
log-exponentially
varying distributions.

Note that partial robustness exists if the conflicting values have a
significant but limited influence on the posterior distribution, as the
conflict grows infinitely.
In contrast, whole robustness is achieved if the influence of the
conflicting values on the posterior distribution
gradually decreases to nothing. To illustrate this, consider
the estimation of a location parameter for a Laplace model (with a
prior of 1). Hence, the posterior mode (or the maximum likelihood estimator)
is the sample median. If, for instance, the sample is $(10, 20, 30, 40,
50, x, x, x, x)$, and we let $x\rightarrow\infty$,
then a wholly robust estimator of the location would be around 30 (the
center of the nonoutlying observations),
while the partially robust sample median estimates the location by 50,
that is, the maximum of the nonoutliers.

This paper goes a step beyond the literature in that it considers
robustness for both location and scale parameters in the same model.
The only other paper that considers Bayesian robustness in a
location--scale model is Andrade and O'Hagan~\cite{r3}.
The essential difference is that partial robustness to a single outlier
is achieved in their paper,
while whole robustness to multiple outliers for both location and scale
estimation is obtained in this paper.

Another distinctive aspect of this paper is the possibility of using
the results of robustness in both frequentist and Bayesian approaches.
Although the model allows us to add prior information on the location
and scale through a very general joint prior density
$\pi(\mu,\sigma)$ [essentially, we only require that $\sigma\pi
(\mu,\sigma)$ is bounded], it is also
possible to choose a noninformative prior such that $\pi(\mu,\sigma
)\propto1/\sigma$. The location and scale parameters can therefore be
estimated in a robust way using either the Bayesian approach
or a frequentist method like maximum likelihood estimation.

This paper is organized as follows. In Section~\ref
{sec-log-regularly}, we introduce
the class of log-regularly varying functions because tail behavior
plays a crucial role in the search of robustness.
Essentially, this class includes functions with a right tail that
exhibits a logarithmic decay, which can be considered a super heavy tail.
As a result, we also define the family of log-regularly varying distributions.

The model with its assumptions is described in Section~\ref
{sec-model}, and
the resolution of conflicts is addressed through the main results of
this paper in Section~\ref{sec-conflict}.
Two simple conditions of robustness are given.
Modeling with a log-regularly varying distribution is the first.
In the second condition,
the number of nonoutlying observations must be larger than the maximum
between the number of small and large outliers.
Results of robustness are asymptotic, where the outlying observations
tend to $-\infty$ or $+\infty$.
Note that the asymptotic nature is about the outliers and not the
sample size, as is usually understood.
Whole robustness is expressed through different types of convergence of
quantities, based on the complete sample,
to quantities based only on the nonoutlying observations, resulting in
a complete rejection of outliers.
We obtain the uniform convergence of the posterior densities, the
convergence in $L_1$, the convergence in distribution
and the uniform convergence of the likelihoods.

In Section~\ref{sec-distribution}, we introduce the family of
log-Pareto-tailed symmetric distributions
that belongs to the larger family of log-regularly varying distributions.
It consists essentially of a symmetric density, such as the standard
normal, with extremities replaced by log-Pareto tails, that is, with
logarithmic decay.
In the presence of outlying observations, the log-Pareto tails ensure
robust inference. Otherwise, the estimation is practically unaffected
by the tails and is determined mostly by the chosen symmetric density.

In Section~\ref{sec-example}, we show that even if the results are
asymptotic, they are still useful in practice with data.
We first illustrate the threshold feature in Section~\ref{sec-threshold}.
When an observation moves away from the nonconflicting values, its
influence on the inference first increases gradually
up to a certain threshold. The conflict then begins, and the model
resolves it by progressively reducing the influence of the
moving observation (now an outlier) to nothing. This built-in feature
is attractive in practice
in that conflict is managed in a sensitive and automatic way.
In Section~\ref{sec-performance}, concurrent estimators are compared
under different scenarios through simulations of observations
to find how they perform in the presence---or absence---of outlying
observations. Nonrobust, partially and wholly robust modeling is considered.
We conclude in Section~\ref{sec-conclusion}, and some proofs are given
in Section~\ref{sec-proof}.

\section{Log-regularly varying functions}\label{sec-log-regularly}

As mentioned in the \hyperref[sec-intro]{Introduction}, tail behavior is crucial for robust
modeling. Hence, we
introduce the class of log-regularly varying functions, as defined in
Desgagn\'e \cite{r8}, following the idea of \mbox{regularly} varying
functions developed by Karamata \cite{r11}.
For each function in Section~\ref{sec-log-regularly}, say $g$, we assume
that $g(z)$ is continuous and strictly positive for $z$ larger than or
equal to a certain constant.

%
\begin{defi}[(Log-regularly varying function)]\label{def-log-regularly}
We say that a measurable function $g$ is \textit{log-regularly
varying} at $\infty$ with
index $\rho\in\re$, written $g\in L_{\rho}(\infty)$, if
\begin{eqnarray*}
&& \forall\varepsilon>0, \forall\tau\ge1,\mbox{ there exists a constant } A(
\varepsilon,\tau)>0\mbox{ such that}
\\[-4pt]
&& z\ge A(\varepsilon,\tau)\mbox{ and }1/\tau\le\nu\le\tau \Rightarrow
\llvert\nu^{\rho}g(z^\nu)/g(z)-1 \rrvert<
\varepsilon.
\end{eqnarray*}
If $\rho=0$, $g$ is said to be \textit{log-slowly varying} at $\infty$.
\end{defi}

In other words, $g\in L_{\rho}(\infty)$ if $g(z^\nu)/g(z)$ converges
to $\nu^{-\rho}$ uniformly in any set $\nu\in[1/\tau,\tau]$ (for
any $\tau\ge1$)
as $z\rightarrow\infty$. The pointwise convergence for any $\nu>0$ follows.

Note that if we define the function $h(z)=g(\ee^z)$, or equivalently
$g(z)=h(\log z)$, we have $g\in L_{\rho}(\infty)$ if and only if $h$
is regularly varying at $\infty$ with index $-\rho$, because
$\lim_{z\rightarrow\infty}h(\nu z)/h(z)=\nu^{-\rho}$. Therefore,
we can obtain different results directly from the theory of regularly
varying functions.
For instance, the functions $\log(\log z)$
and 1 are both log-slowly varying at $\infty$ since $\log z$ and 1 are
slowly varying.

%
\begin{prop}[(Equivalence)]\label{prop-equivalence}
For any $\rho\in\re$, we have $g\in L_{\rho}(\infty)$ if and only if
there exists a constant $A>1$ and a function $s\in L_{0}(\infty)$ such
that for
$z\ge A$, $g$ can be written as
\[
g(z)=(\log z)^{-\rho} s(z).
\]
\end{prop}

\begin{pf}It is well known that if a function $h$ is regularly varying
at $\infty$ with index $-\rho$, it can be represented as
$h(z)=z^{-\rho}l(z)$, where $l$ is some slowly varying function. It is
equivalent to say that $g\in L_{\rho}(\infty)$, where
\[
g(z)=h(\log z)=(\log z)^{-\rho}l(\log z)=(\log z)^{-\rho}s(z),
\]
with $s(z)=l(\log z)\in L_{0}(\infty)$.
\end{pf}

The next proposition establishes the asymptotic dominance of
a logarithmic function over a log-slowly varying function.

%
\begin{prop}[(Dominance)]\label{prop-dominance}
If $s\in L_{0}(\infty)$ and $g\in L_{\rho}(\infty)$, then for all
$\delta>0$, there exists a constant $A(\delta)> 1$ such that
$z\ge A(\delta)\Rightarrow$
\[
(\log z)^{-\delta}< s(z) < (\log z)^{\delta}\quad\mbox{and}\quad(
\log z)^{-\rho-\delta}< g(z) < (\log z)^{-\rho+\delta}.
\]
\end{prop}

\begin{pf}
It is well known that if $l$ is slowly varying, then for every $\delta
>0$, we have $z^{-\delta}l(z)\rightarrow0$ and
$z^{\delta}l(z)\rightarrow\infty$ as $z\rightarrow\infty$. It
follows that $z^{-\delta}<l(z)<z^{\delta}$ for $z$ sufficiently large.
If we replace $z$ by $\log z$ and we set $s(z)=l(\log z)$, then $s\in
L_0(\infty)$, and we obtain that $(\log z)^{-\delta}s(z)\rightarrow
0$ and
$(\log z)^{\delta}s(z)\rightarrow\infty$ as $\log z\rightarrow
\infty$ (or equivalently $z\rightarrow\infty$) and $(\log
z)^{-\delta}<s(z)<(\log z)^{\delta}$ for $z$ sufficiently large.
Since we can write $g(z)=(\log z)^{-\rho} s(z)$, the second part of
the proposition follows directly.
\end{pf}

The index $\rho$ can be interpreted as a measure of the tail's
thickness or as a tail index, which is useful for the ordering of
different tails.
The function with the smallest tail index $\rho$ has the heaviest
tail. More formally,
we can verify that if $g_1\in L_{\rho_1}(\infty)$ and $g_2\in L_{\rho
_2}(\infty)$, then
$\rho_1>\rho_2\Rightarrow g_1(z)/g_2(z)\rightarrow0$ as
$z\rightarrow\infty$.
The tail index $\rho$ is also useful to determine if $(1/z)g(z)$ is
integrable, where $g(z)\in L_{\rho}(\infty)$,
as described in the next proposition.

%
\begin{prop}[(Integrability)]\label{prop-L-integrability}
If $g(z)\in L_{\rho}(\infty)$,
then there exists a constant $A>0$ such that
$(1/z)g(z)$ is integrable on $z\ge A$, if and only if:

\begin{longlist}[(ii)]
\item[(i)] $\rho>1$,

\item[(ii)] $\rho=1$, with the log-slowly varying part of $g(z)$
having a sufficiently fast decay [e.g., faster than $(\log(\log
z))^{-\beta}$, with $\beta>1$].
\end{longlist}
\end{prop}
\begin{pf}
If we define $h$ such that $g(z)=h(\log z)$, and we choose $A$
sufficiently large, then $h$ is regularly varying at $\infty$ with
index $-\rho$, and we have
\[
\int_{A}^\infty(1/z)g(z) \,dz =\int
_{A}^\infty(1/z)h(\log z) \,dz=\int
_{\log A}^\infty h(u) \,du=\int_{\log A}^\infty
u^{-\rho}l(u) \,du,
\]
where $l$ is slowly varying. For any $\delta>0$, if $A$ is
sufficiently large, we have $u^{-\delta}<l(u)<u^{\delta}$. Therefore,
the integral
exists if $\rho>1$ and does not if $\rho<1$.

If $\rho=1$, we see that the decay of $l$ determines the existence of
the integral.
If, for instance, $l(u)<(\log u)^{-\beta}$ or $s(z)=l(\log u)<(\log
(\log u))^{-\beta}$, with $\beta>1$ and $s\in L_0(\infty)$, then the
integral exists.
Instead, if $l(u)>(\log u)^{-\beta}$ or $s(z)=l(\log u)>(\log(\log
u))^{-\beta}$, with $\beta<1$ and $s\in L_0(\infty)$, then the
integral does not exist.
\end{pf}

In particular, if $f$ is a continuous symmetric probability density
function defined on $\re$ such that
$g(z)=z f(z)\in L_{\rho}(\infty)$, we know from
Proposition~\ref{prop-L-integrability} that a tail index $\rho>1$
is sufficient to guarantee that $f$ is proper and that $\rho\ge1$ is
a necessary condition.
This leads us to the next definition.

%
\begin{defi}[(Log-regularly varying distribution)]\label
{def-log-regularly-distribution}
A random variable $Z$ and its distribution are said to be log-regularly
varying with index $\rho\ge1$
if their symmetric density $f$ is such that $z f(z)\in L_{\rho}(\infty)$.
\end{defi}

Using Propositions~\ref{prop-equivalence} and \ref{prop-dominance},
this means that for all $\delta>0$ and $\llvert z\rrvert$
larger than a certain constant,
the symmetric (with respect to 0) density $f$ of a log-regularly
varying distribution with index $\rho$ can be written as
$f(z)=(1/\llvert z\rrvert)(\log\llvert z\rrvert)^{-\rho
} s(\llvert z\rrvert)$, where $s\in L_{0}(\infty)$
can be bounded by $(\log\llvert z\rrvert)^{-\delta}$ and $(\log
\llvert z\rrvert)^{\delta}$.
Such a density with logarithmic decaying tails can be referred to as a
super heavy-tailed distribution.

In the next proposition, we see the asymptotic impact of a
location--scale transformation on a log-regularly varying function $g$
and the density $f$ of a log-regularly varying distribution.
Mostly, it is another way to express tail thickness.

%
\begin{prop}[(Location--scale transformation)]\label
{prop-location--scale-transformation}
If $g(z)=z f(z)\in L_{\rho}(\infty)$, then we have, as $z\rightarrow
\infty$,
\[
g \bigl((z-\mu)/\sigma \bigr)/g(z)\rightarrow1\quad\mbox{and}\quad (1/\sigma) f
\bigl((z-\mu)/\sigma \bigr)/f(z)\rightarrow1,
\]
uniformly on $(\mu,\sigma)\in[-\lambda,\lambda]\times[1/\tau
,\tau]$, for any $\lambda\ge0$ and $\tau\ge1$.
\end{prop}
\begin{pf}
Using $g\in L_{\rho}(\infty)$ with Proposition~\ref
{prop-equivalence}, there exists a
function $s\in L_{0}(\infty)$ such that $g(z)=(\log z)^{-\rho} s(z)$,
if $z$ is large enough.
Therefore, for any chosen $\lambda\ge0$ and $\tau\ge1$, if $z$ is
sufficiently large, we have
\begin{eqnarray*}
\frac{g((z-\mu)/\sigma)}{g(z)}& =& \biggl(\frac{\log((z-\mu
)/\sigma)}{\log
z} \biggr)^{-\rho}
\frac{s((z-\mu)/\sigma)}{s(z)}.
\end{eqnarray*}
It is purely algebraic to show that the term $(\log((z-\mu)/\sigma
))/(\log z)$ converges to~1 uniformly on any set
$(\mu,\sigma)\in[-\lambda,\lambda]\times[1/\tau, \tau]$ as
$z\rightarrow\infty$.

Finally, we want to show that $s((z-\mu)/\sigma)/s(z)$ converges to 1
uniformly on any set
$(\mu,\sigma)\in[-\lambda,\lambda]\times[1/\tau, \tau]$ as
$z\rightarrow\infty$, or equivalently that $s(y)/s(z)$ converges to 1
uniformly on $y\in[(z-\lambda)/\tau,(z+\lambda)\tau]$.
We observe that for any chosen $\lambda\ge0$ and $\tau\ge1$, if $z$
is sufficiently large, we have
\[
z^{1/2}\le(z-\lambda)/\tau\le(z+\lambda)\tau\le z^2.
\]
Therefore, it suffices to show that $s(y)/s(z)$ converges to 1
uniformly on $y\in[z^{1/2},z^2]$, or equivalently, that
$s(z^\nu)/s(z)$ converges to 1 uniformly on any set $\nu\in[1/2,2]$,
which is the case since $s\in L_{0}(\infty)$.
The second part of the proposition follows directly.
\end{pf}

\section{Resolution of conflicts in a location--scale parameter
model}\label{sec-robustness}

\subsection{Model}\label{sec-model}
\begin{longlist}[(ii)]
\item[(i)]Let $X_1,\ldots,X_n$ be $n$ random variables conditionally
independent given $\mu$ and $\sigma$ with their conditional densities
given by
\[
X_i \mid\mu,\sigma\simdist(1/\sigma)f \bigl((x_i-\mu)/
\sigma \bigr);
\]
\item[(ii)]the joint prior density of $\mu$ and $\sigma$ is given by
$\mu,\sigma\simdist\pi(\mu,\sigma)$,
where $n\ge2, x_1,\ldots,x_n,\mu\in\re$, $\sigma>0$.
\end{longlist}

We assume that the prior $\pi(\mu,\sigma)$ is nonnegative on $\re
$, and
the only other required assumption is that $\sigma\pi(\mu,\sigma)$
is bounded. Note that in particular,
if we have no prior information or if we use the model in a frequentist
approach, then we set
$\pi(\mu,\sigma)\propto1/\sigma$, an improper joint prior density
which can be considered as noninformative.

We assume that $f$ is a proper density that is continuous and strictly
positive on $\re$.
In addition, we assume it is symmetric with respect to the origin.
We also assume that both tails of $\llvert z\rrvert f(z)$ are
monotonic, which means
that the tails of $f(z)$ are also monotonic.
Note that monotonicity of the tails of $f(z)$ and $\llvert z\rrvert
f(z)$ means that
there exists
a constant $M\ge0$ such that
%
%
\begin{equation}
\label{eqn-monotonic} \llvert y\rrvert\ge\llvert z\rrvert\ge M\qquad\mbox{implies that
}f(y)\le f(z)\quad\mbox{and}\quad\llvert y\rrvert f(y)\le \llvert z\rrvert f(z).
\end{equation}
It follows that $f(z)$ and $\llvert z\rrvert f(z)$
are bounded on the real line, with a limit of 0 in their tails as
$\llvert z\rrvert\rightarrow\infty$.
Hence, considering also the prior, we can define the constant~$B$ as follows:
%
%
\begin{equation}
\label{eqn-B} B=\max \Bigl\{\sup_{z\in\re}f(z),\sup
_{z\in\re}\llvert z\rrvert f(z),\sup_{\mu\in\re,\sigma>0}\sigma
\pi(\mu,\sigma) \Bigr\}.
\end{equation}
These conditions are referred to below as the conditions of regularity
on $f$.
The density $f$ can possess other parameters than location and scale,
such as a shape parameter, but they are assumed to be known.

We study robustness of the estimation of $\mu$ and $\sigma$ in the
presence of outliers.
The nature of the results is asymptotic, in the sense that some $x_i$ are
going to $-\infty$ or $+\infty$. We want to find sufficient
conditions to obtain whole robustness, that is, a
complete rejection of the outliers.

Among the $n$ observations, denoted by $\mathbf{x_n}=(x_1,\ldots,x_n)$,
we assume that $k\ge2$ of them, denoted by the vector $\mathbf{x_k}$,
form a group of nonoutlying observations.
We assume that $l$ of them
are considered as left outliers (smaller than
the nonoutliers) and $r$ of them
are considered as right outliers (larger than the nonoutliers), with $k+l+r=n$.

For $i=1,\ldots,n$, we define three binary functions $k_i, l_i$ and
$r_i$ as follows.
If $x_i$ is a nonoutlying observation, we set $k_i=1$; if it is a left
outlier, we set $l_i=1$; and if it is a right outlier, we set $r_i=1$.
These functions are set to 0 otherwise.
We have $k_i+l_i+r_i=1$ for $i=1,\ldots,n$, with $\sum_{i=1}^n
k_i=k$, $\sum_{i=1}^n l_i=l$ and $\sum_{i=1}^n r_i=r$.

We assume that each outlier is going to $-\infty$ or $+\infty$ at its
own specific rate, to the extent that the ratio of two outliers is bounded.
We can write
\[
x_i=a_i+b_i \omega,
\]
for $i=1,\ldots,n$, where $a_i$ and $b_i$ are some constants such that
$a_i\in\re$ and:
\begin{longlist}[(iii)]
\item[(i)] $b_i=0$ if $k_i=1$;

\item[(ii)] $b_i<0$ if $l_i=1$;

\item[(iii)] $b_i>0$ if $r_i=1$;
\end{longlist}
and we let $\omega\rightarrow\infty$. Note that if multiple outliers
share the same $b_i$, they move as a block at the same rate.

Let the joint posterior density of $\mu$ and $\sigma$ be denoted by
$\pi(\mu,\sigma\mid\mathbf{x_n})$
and the marginal density of $X_1,\ldots,X_n$ be denoted by $m(\mathbf
{x_n})$, with
\[
\pi(\mu,\sigma\mid\mathbf{x_n})= \bigl[m(\mathbf{x_n})
\bigr]^{-1}\pi(\mu,\sigma)\prod_{i=1}^n
(1/\sigma)f \bigl((x_i-\mu)/\sigma \bigr).
\]
Let the joint posterior density of $\mu$ and $\sigma$ considering
only the nonoutlying observations $\mathbf{x_k}$
be denoted by $\pi(\mu,\sigma\mid\mathbf{x_k})$
and its corresponding marginal density be denoted by $m(\mathbf
{x_k})$, with
\[
\label{eqn-non-outlier} \pi(\mu,\sigma\mid\mathbf{x_k})= \bigl[m(
\mathbf{x_k}) \bigr]^{-1}\pi(\mu,\sigma)\prod
_{i=1}^n \bigl[(1/\sigma)f \bigl((x_i-
\mu)/\sigma \bigr) \bigr]^{k_i}.
\]

The likelihood functions can be found by setting $\pi(\mu,\sigma
)\propto1/\sigma$ and letting
\[
\likely(\mu,\sigma\mid\mathbf{x_n})\propto\sigma\pi(\mu,\sigma \mid
\mathbf{x_n})\quad\mbox{and}\quad\likely(\mu,\sigma\mid
\mathbf{x_k})\propto\sigma\pi(\mu,\sigma\mid\mathbf{x_k}).
\]

%
\begin{prop}\label{proposition-proper}
Considering the Bayesian context given in Section~\ref{sec-model},
the joint posterior densities $\pi(\mu,\sigma\mid\mathbf{x_k})$
and $\pi(\mu,\sigma\mid\mathbf{x_n})$ are proper.
\end{prop}

The proof of Proposition~\ref{proposition-proper} is given in
Section~\ref{sec-proof}.

\subsection{Resolution of conflicts}\label{sec-conflict}

The results of robustness are now given.

%
\begin{teo}\label{teo-main}
Consider the model and context described in Section~\ref{sec-model},
and assume that the conditions of regularity on $f$ are satisfied.
If we have:
\begin{longlist}[(ii)]
\item[(i)]$z f(z)\in L_{\rho}(\infty)$ [$z f(z)$ is
log-regularly varying at $\infty$ with index $\rho\ge1$],
\item[(ii)]$k>\max(l,r)$,
\end{longlist}
then we obtain the following results:
\begin{longlist}[(a)]
\item[(a)]
\[
\lim_{\omega\rightarrow\infty}\frac{m(\mathbf{x_n})}{\prod_{i=1}^{n}[f(x_i)]^{l_i+r_i}}= m(\mathbf{x_k}).
\]
\item[(b)]
\[
\lim_{\omega\rightarrow\infty}\pi(\mu,\sigma\mid\mathbf {x_n})=\pi(
\mu,\sigma\mid\mathbf{x_k}),
\]
uniformly on $(\mu,\sigma)\in[-\lambda,\lambda]\times[1/\tau
,\tau]$, for any $\lambda\ge0$ and $\tau\ge1$.\vspace*{1pt}
\item[(c)]
\[
\lim_{\omega\rightarrow\infty}\int_{0}^{\infty}\int
_{-\infty
}^{\infty} \bigl\llvert\pi(\mu,\sigma\mid
\mathbf{x_n})-\pi(\mu,\sigma\mid\mathbf{x_k}) \bigr\rrvert
\,d\mu\,d\sigma= 0.
\]
\item[(d)]As $\omega\rightarrow\infty$,
\[
\mu,\sigma\mid\mathbf{x_n} \convdist\mu,\sigma\mid
\mathbf{x_k},
\]
and in particular
\[
\mu\mid\mathbf{x_n} \convdist\mu\mid\mathbf{x_k}\quad
\mbox{and}\quad\sigma\mid\mathbf{x_n} \convdist\sigma\mid
\mathbf{x_k}.
\]
\item[(e)]
\[
\lim_{\omega\rightarrow\infty}\likely(\mu,\sigma\mid\mathbf{x_n})=
\likely(\mu,\sigma\mid\mathbf{x_k}),
\]
uniformly on $(\mu,\sigma)\in[-\lambda,\lambda]\times[1/\tau
,\tau]$, for any $\lambda\ge0$ and $\tau\ge1$.
\end{longlist}
\end{teo}

Proof of result (a) is substantial and therefore is given in
Section~\ref{sec-proof}. This is, however, the crucial part in the
proof of Theorem~\ref{teo-main}.
\begin{pf*}{Proof of result \textup{(b)}}
Consider $(\mu,\sigma)$ such that $\pi(\mu,\sigma)>0$ [the proof
for the case $(\mu,\sigma)$ such that $\pi(\mu,\sigma)=0$ is trivial].
We have, as $\omega\rightarrow\infty$,
\begin{eqnarray*}
\frac{\pi(\mu,\sigma\mid\mathbf{x_n})}{\pi(\mu,\sigma\mid
\mathbf{x_k})} &=& \frac{m(\mathbf{x_k})}{m(\mathbf{x_n})} \frac
{\pi(\mu,\sigma)\prod_{i=1}^{n}(1/\sigma)f((x_i-\mu)/\sigma)}{
\pi(\mu,\sigma)\prod_{i=1}^{n} [(1/\sigma)f((x_i-\mu
)/\sigma) ]^{k_i}}
\\
&=& \frac{m(\mathbf{x_k})}{m(\mathbf{x_n})}\prod_{i=1}^{n}
\bigl[(1/\sigma)f \bigl((x_i-\mu)/\sigma \bigr) \bigr]^{l_i+r_i}
\\
&=& \frac{m(\mathbf{x_k})\prod_{i=1}^{n}[f(x_i)]^{l_i+r_i}}{m(\mathbf
{x_n})}\cdot\prod_{i=1}^{n}
\biggl[\frac{(1/\sigma)f((x_i-\mu
)/\sigma)}{
f(x_i)} \biggr]^{l_i+r_i}\rightarrow1.
\end{eqnarray*}
The first part of the last term does not depend on $\mu$ and $\sigma$
and converges to 1 as $\omega\rightarrow\infty$, using result (a).
The second part of the last term also converges to~1 uniformly
in any set $(\mu,\sigma)\in[-\lambda,\lambda]\times[1/\tau, \tau
]$ using Proposition~\ref{prop-location--scale-transformation}.
Furthermore, since $f$ and $\sigma\pi(\mu,\sigma)$ are bounded,
$\pi(\mu,\sigma\mid\mathbf{x_k})$ is also bounded on any set $(\mu
,\sigma)\in[-\lambda,\lambda]\times[1/\tau, \tau]$.
Then we have
\begin{eqnarray*}
&& \bigl\llvert\pi(\mu,\sigma\mid\mathbf{x_n})-\pi(\mu,\sigma \mid
\mathbf{x_k}) \bigr\rrvert
\\
&&\qquad =\pi(\mu,\sigma\mid\mathbf{x_k})
\biggl\llvert \frac {\pi(\mu,\sigma\mid\mathbf{x_n})}{\pi(\mu,\sigma
\mid\mathbf {x_k})}-1\biggr\rrvert\rightarrow0 \qquad\mbox{as }
\omega \rightarrow\infty.
\end{eqnarray*}\upqed
\end{pf*}
\begin{pf*}{Proofs of results \textup{(c)} and \textup{(d)}}
We can use Scheff\'{e}'s theorem \cite{r17} directly to prove results
(c) and (d). Using Proposition~\ref{proposition-proper}, we know that
$\pi(\mu,\sigma\mid\mathbf{x_k})$ and $\pi(\mu,\sigma\mid
\mathbf{x_n})$ are proper.
Using result (b), we have that $\pi(\mu,\sigma\mid\mathbf
{x_n})\rightarrow\pi(\mu,\sigma\mid\mathbf{x_k})$
pointwise as $\omega\rightarrow\infty$ for any $\mu\in\re$ and
$\sigma>0$, as a result of the uniform convergence. The conditions of
Scheff\'{e}'s theorem are then satisfied, and
we obtain the convergence in $L_1$ given in result (c) as well as the
following result:
\[
\lim_{\omega\rightarrow\infty}\int_{E}\pi(\mu,\sigma\mid
\mathbf{x_n}) \,d\mu\,d\sigma= \int_{E}\pi(\mu,
\sigma\mid\mathbf{x_k}) \,d\mu\,d\sigma,
\]
uniformly for all rectangles $E$ in $\re\times\re^{+}$.
\end{pf*}
\begin{pf*}{Proof of result \textup{(e)}}
It suffices to write the likelihood functions as $\likely(\mu,\sigma
\mid\mathbf{x_n})\propto\sigma\pi(\mu,\sigma\mid\mathbf{x_n})$ and
$\likely(\mu,\sigma\mid\mathbf{x_k})\propto\sigma\pi(\mu
,\sigma\mid\mathbf{x_k})$ with $\pi(\mu,\sigma)\propto1/\sigma
$, and result (e) follows directly from result (b).
\end{pf*}

An attractive feature of Theorem~\ref{teo-main} is the simplicity of
its only two sufficient conditions. Condition (i) says that
modeling must be done using density $f$ of a log-regularly varying
distribution with index $\rho\ge1$; see Definition~\ref
{def-log-regularly-distribution}. Note that it involves only the tails
of the function $\llvert z\rrvert f(z)$.
Essentially, the decay of the tails must be logarithmic.
For that purpose, in the next section we introduce the family of
log-Pareto-tailed symmetric distributions
that belong to the family of log-regularly varying distributions.

Condition (ii) requires that $k>l$ and $k>r$. For instance, a group of
$k=6$ nonoutlying observations is sufficient
to ensure the rejection of $l=5$ outliers at left and $r=5$ at right.
The nonoutlying group must be the most important, which is rather intuitive.
The most demanding case occurs when all outliers are on the same side
(e.g., $l=0$). Condition (ii) can then be written as $k>n/2$,
which means that the nonoutliers must represent more than half of the sample.
A few numerical simulations tend to confirm our expectation that a
larger difference between $k$
and $\max(l,r)$ results in a faster rejection of the outliers.

The asymptotic behavior of the marginal $m(\mathbf{x_n})$ is given in
result (a). This fundamental result
is probably of more theoretical than practical interest because it
leads to results (b) to (e).
The asymptotic behavior of the posterior density is given in result (b).
The posterior considering the entire sample converges to the posterior
considering only the
$k$ nonoutlying observations, uniformly in any set $(\mu,\sigma)\in
[-\lambda,\lambda]\times[1/\tau,\tau]$.
The outliers are then completely rejected as they are going to plus or
minus infinity.
We also obtain the pointwise convergence.

In result (c), we obtain the convergence in $L_1$ of the posterior
densities considering the entire sample
to the posterior considering only the nonoutlying observations. In
result (d), we obtain the convergence in distribution, that is
$\Pr(\mu,\sigma\in E \mid\mathbf{x_n})$ converges to $\Pr(\mu
,\sigma\in E \mid\mathbf{x_k})$ as $\omega\rightarrow\infty$,
uniformly for all rectangles $E$ in $\re\times\re^{+}$.
Because the convergence is uniform, this is actually a stronger result
than the convergence in distribution, which requires only pointwise convergence.
We also obtain the convergence in distribution of the posterior
marginal distributions.
Therefore, any estimation of $\mu$ and $\sigma$ based on posterior
quantiles or Bayesian credible intervals is robust to outliers.

In result (e), the likelihood considering the entire sample converges
to the likelihood considering only the
nonoutlying observations, uniformly in any set $(\mu,\sigma)\in
[-\lambda,\lambda]\times[1/\tau,\tau]$.
It follows that the maximum of $\likely(\mu,\sigma\mid\mathbf
{x_n})$ converges to the maximum of
$\likely(\mu,\sigma\mid\mathbf{x_k})$, and therefore the maximum
likelihood estimates also converge, as $\omega\rightarrow\infty$.

\section{The family of log-Pareto-tailed symmetric
distributions}\label{sec-distribution}

As stated in Theorem~\ref{teo-main}, modeling with a log-regularly
varying distribution is one of the conditions of robustness.
However, such a distribution is super heavy-tailed, and the usual
densities defined on $\re$ are light or heavy-tailed.
Therefore, we introduce in this section the family of log-Pareto-tailed
symmetric distributions
that belongs to the larger family of log-regularly varying
distributions. Given that the conditions of robustness involve only the
tails of density $f(z)$,
the proposed solution consists in altering a symmetric density, such as
the usual normal, uniform or
Student's $t$ distributions, by replacing its extremities with
log-Pareto tails, that is, a function proportional to
$\llvert z\rrvert^{-1}(\log\llvert z\rrvert)^{-\beta}$,
with $\beta>1$. This idea comes from
the generalized exponential power (GEP) distribution, a family
introduced by
Angers \cite{r5} and revisited in more detail by Desgagn\'e and Angers
\cite{r9}.
The GEP density is essentially a uniform density in the center with a
large spectrum of tail behavior, classified in types I to V, from light
to super heavy-tailed.
In particular, the GEP of type~V is a log-regularly varying
distribution because its density has log-Pareto tails. We propose here
to generalize the GEP distribution of type~V
to the family of log-Pareto-tailed symmetric distributions by using any
symmetric densities in the center instead of limiting the choice to the
uniform density.

%
\begin{defi}
A random variable $Z$ has a log-Pareto-tailed symmetric distribution if its
density is given by
\begin{eqnarray*}
&& f(z\mid\bolds\phi, \alpha,\beta)
\\
&&\qquad = K_{(\bolds\phi,\alpha,\beta)} \biggl(g(z\mid\bolds\phi)
\mathbh {1}_{[-\alpha,\alpha]}(z) + g(\alpha\mid\bolds\phi)\frac{\alpha}{\llvert z\rrvert} \biggl(
\frac{\log
\alpha}{\log\llvert z\rrvert} \biggr)^\beta\mathbh{1}_{(\alpha,\infty
)} \bigl(\llvert z
\rrvert \bigr) \biggr),
\end{eqnarray*}
where $z\in\re$, $\alpha> 1$, $\beta>1$, $\mathbh{1}_{A}(\cdot)$
is an indicator function, and $g(\cdot\mid\bolds\phi)$ is any
density that is symmetric with respect to the origin,
continuous and strictly positive on $[-\alpha,\alpha]$,
with its vector of parameters given by $\bolds\phi\in
\bolds{\Phi}$.
The normalizing constant is given by
\[
K_{(\bolds\phi,\alpha,\beta)}=\frac{(\beta-1)}{(2 G(\alpha
\mid\bolds\phi)-1)(\beta-1)+2 g(\alpha\mid\bolds\phi
)\alpha\log\alpha},
\]
where $G(\alpha\mid\bolds\phi)=\int_{-\infty}^{\alpha
}g(u\mid\bolds\phi) \,du$.
\end{defi}

In particular, if $g(z\mid\bolds\phi)$ is a normal density, we
say that the random variable $Z$ has a log-Pareto-tailed normal distribution.
If $g(z\mid\bolds\phi)$ is a Student's $t$ density, we say that
$Z$ has a log-Pareto-tailed Student's $t$ distribution, and so on.
The core of the density $f(z\mid\bolds\phi,\alpha,\beta)$ is
located between $-\alpha$ and $\alpha$, and the tails are positioned
in the area $\llvert z\rrvert> \alpha$.
Tail thickness is controlled with the parameter $\beta$.
This density satisfies the condition of robustness required in
Theorem~\ref{teo-main}, since for $\llvert z\rrvert>\alpha$, we have
\[
\llvert z\rrvert f(z\mid\bolds\phi,\alpha,\beta)\propto \bigl(\log\llvert z
\rrvert \bigr)^{-\beta
}\in L_\beta(\infty).
\]
All conditions of regularity assumed in Section~\ref{sec-model} are
satisfied as well.
The density $f(z\mid\bolds\phi,\alpha,\beta)$ is continuous
and strictly positive on $\re$, proper (see Proposition~\ref
{prop-L-integrability}) and
symmetric with respect to the origin. Furthermore, both tails of
$\llvert z\rrvert
f(z\mid\bolds\phi,\alpha,\beta)$ are monotonic.

In practice, choosing parameters $\alpha$ and $\beta$ directly is not
necessarily an intuitive task. It could be easier to choose other
indirect but related
quantities. Here is an interesting strategy in five steps: a~practitioner first chooses his favorite symmetric density $g(z\mid
\bolds\phi)$
and its vector of parameters $\bolds\phi$ (other than the
location and scale parameters $\mu$ and $\sigma$, which will be added
later), such as the $N(0,1)$.
The second step consists in setting the normalizing constant
$K_{(\bolds\phi,\alpha,\beta)}$ to 1. The desirable
consequence is that
the core (between $-\alpha$ and $\alpha$) of the density $f(z\mid
\bolds\phi,\alpha,\beta)$ becomes exactly the density $g(z\mid
\bolds\phi)$,
the familiar density of the user. The third step consists in choosing
the mass of the core, which is defined as
\[
q=\Pr(-\alpha\le Z \le\alpha\mid\bolds\phi,\alpha,\beta).
\]

For instance, we could choose $q=0.95$, which leaves 2.5\% of the mass
in each tail. Then, the density $f(z\mid\bolds\phi,\alpha,\beta)$
would be exactly the $N(0,1)$ density for 95\% of its mass located in
the center. The following steps are done automatically.
Given that $K_{(\bolds\phi,\alpha,\beta)}$ has been set to 1,
it follows that $q=2 G(\alpha\mid\bolds\phi)-1$.
However, to ensure that $\alpha>1$ as required, we must choose $q>2
G(1\mid\bolds\phi)-1$.
If the last equality is rearranged, it leads us to the fourth step,
which consists in calculating $\alpha$ as follows:
\[
\label{eqn-alpha} \alpha=G^{-1} \biggl(\frac{1+q}{2} \biggm| \bolds
\phi \biggr).
\]
For example, a $N(0,1)$ with $q=0.95$ generates a value of $\alpha=1.96$.
Finally, we calculate $\beta$ in the fifth step as follows:
\[
\label{eqn-beta}\beta= 1+\frac{2 g(\alpha\mid\bolds\phi) \alpha
\log\alpha}{1-q}.
\]
Note that this equation is consistent with a normalizing constant of 1,
and it satisfies $\beta>1$ since $\alpha>1$. Our example gives a
value of $\beta=4.08$.

%
\begin{figure}

\includegraphics{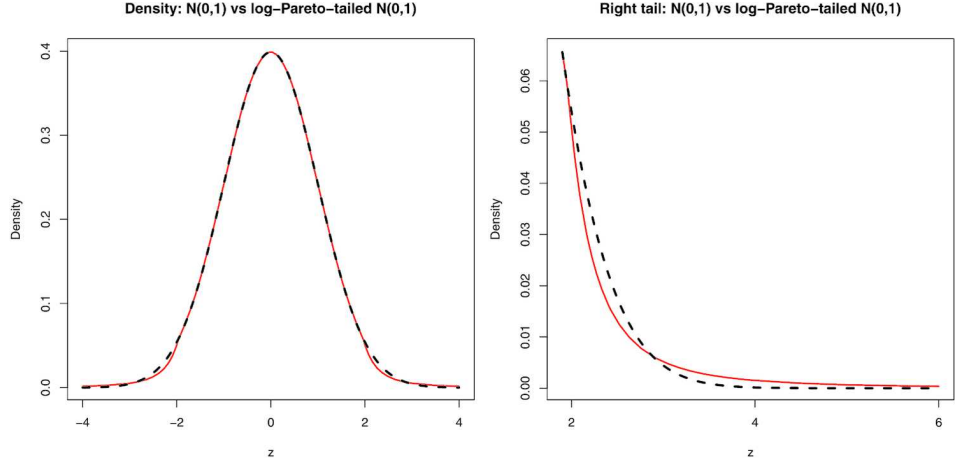}

\caption{A comparison between the standard normal (dashed line)
and log-Pareto-tailed standard normal (solid line)
densities.}\label{figlog-Pareto-density}
\end{figure}

We compare in Figure~\ref{figlog-Pareto-density} the standard normal
density (dashed line) to a log-Pareto-tailed standard normal
density (solid line), with
$q=0.95$,\break $K_{(\bolds\phi,\alpha,\beta)}=1$, $\alpha=1.96$
and $\beta=4.08$. Both densities are identical between $-\alpha$ and
$\alpha$,
but differ in the tails.

Simulation of observations from a log-Pareto-tailed symmetric
distribution is easy using the inverse transformation method.
It is described in detail in Section~3.4 of Desgagn\'e and Angers \cite
{r9} for the log-Pareto-tailed uniform distribution
(labeled GEP density of type~V in their paper). It is straightforward
to generalize it to other symmetric densities $g(\cdot\mid\bolds
\phi)$.

Of course, we can add location and scale parameters, denoted,
respectively, by $\mu\in\re$ and $\sigma>0$, to the density
$f(z\mid
\bolds\phi,\alpha,\beta)$. We obtain
\begin{eqnarray*}
&& (1/\sigma)f \bigl((z-\mu)/\sigma\mid\bolds\phi,\alpha,\beta \bigr)
\\
&&\qquad =\cases{
K_{(\bolds\phi,\alpha,\beta)}(1/\sigma) g \bigl((z-\mu )/\sigma\mid \bolds\phi \bigr),\qquad\mbox{if $\displaystyle\mu- \alpha\sigma\le z\le\mu+ \alpha\sigma$},
\vspace*{3pt}\cr
\displaystyle K_{(\bolds\phi,\alpha,\beta)}g(\alpha\mid
\bolds\phi)\frac
{\alpha}{\llvert z-\mu\rrvert} \biggl( \frac{\log\alpha}{\log
(\llvert z-\mu
\rrvert/\sigma)} \biggr)^\beta,
\cr
\hspace*{154pt}\quad \mbox{if $\llvert z-\mu\rrvert \ge \alpha\sigma$}.}
\end{eqnarray*}
\label{densitylogpareto}

Note that when this density is used in the context of robustness
described in Section~\ref{sec-conflict}, the parameters
$\bolds\phi$, $\alpha$ and $\beta$ are assumed to be known.
The inference is done on the location and scale parameters only.

\section{Example}\label{sec-example}

In this section,
the asymptotic results of robustness found in Theorem~\ref{teo-main}
are confronted with data. Without loss of generality,
we choose the improper and
noninformative joint prior density $\pi(\mu,\sigma)\propto1/\sigma
$. Hence, both the Bayesian and frequentist approaches can be used.

We first
illustrate in Section~\ref{sec-threshold} the behavior of different
estimators of the location and scale parameters when one observation
moves from 0 to 100, given that
the rest of data lie between $-$10 and 10. For the estimator based on
robust modeling provided by Theorem~\ref{teo-main}, we observe an interesting
feature that we call the threshold. The influence of the moving
observation on the inference increases until a certain threshold. Then
the nature of this observation
gradually changes to become more and more outlying, as its influence
decreases and eventually completely disappears. In Section~\ref
{sec-performance},
the performances of concurrent estimators are compared for different
scenarios. We consider simulation of observations
from the normal as well as from contaminated normal distributions, to
see how the estimators perform in the presence---or absence---of
outliers. The mean square error
is calculated as the measure of performance.

\subsection{Illustration of the threshold}\label{sec-threshold}

We consider a sample of size $n=22$ given by $\mathbf{x_n}=(\mathbf
{x_k},\omega)$, where the $k=21$ nonoutlying observations are represented
by $\mathbf{x_k}=(-10,-9, \ldots,-1,0,1,\ldots,9,10)$.
We study the impact of moving the observation $\omega$ from 0 to 100
on the location--scale parameter inference based on the maximum
likelihood estimator (MLE)
calculated for three different densities $f$, in accordance with the
model described in Section~\ref{sec-model}.
Note that results using the Bayesian marginal posterior median are very similar.
Naturally, the standard normal density has been chosen as the
nonrobust model. The corresponding MLE
are then the usual sample mean and (biased) sample standard deviation.

The log-Pareto-tailed standard normal density, as illustrated in
Figure~\ref{figlog-Pareto-density},
is also studied. We have chosen $q=0.95$, $\alpha=1.96$ and $\beta
=4.08$, as discussed in Section~\ref{sec-distribution}.
This modeling leads to complete rejection of the outlier, as described
by Theorem~\ref{teo-main}.
We also examined other values of $q$ (the values of $\alpha$ and
$\beta$ are calculated automatically using the proposed algorithm in
Section~\ref{sec-distribution}).
If we choose a larger value of $q$, then the density is closer to the
$N(0,1)$, and the same goes for the inference in the absence of outliers.
However, the threshold of robustness increases. The choice of 0.95
appeared to be well balanced for good inference with and without outliers.

%
\begin{figure}

\includegraphics{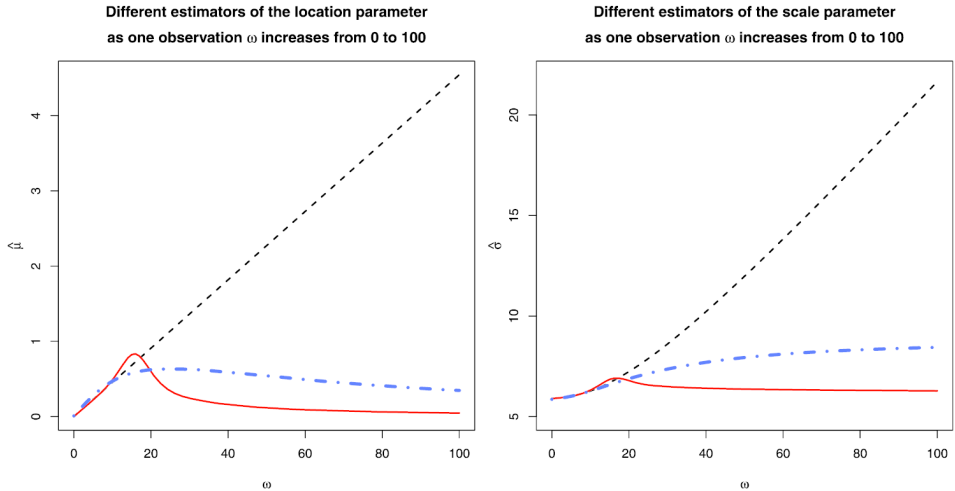}

\caption{Estimation of the location (left graph) and scale (right
graph) parameters for the normal model (dashed lines), the
log-Pareto-tailed normal model (solid lines) and the Student's $t$
model (dotted-dashed lines), using the MLE.}\label{figrobustnesslocation}
\end{figure}

The third density $f$ considered is the Student's $t$, a common choice
for robust modeling. This density satisfies the conditions of robustness
given in Andrade and O'Hagan \cite{r3} (which lead to partial
robustness concerning the scale parameter), but not the conditions of
whole robustness given in Theorem~\ref{teo-main}. The degrees of
freedom has been set to 10, again to search for balance between good
inference with and without outliers.
An implicit scale parameter of 0.964 (other than $\sigma$) has been
added to match its interquartile range to that of the two densities
considered above.

Robustness for the three models is illustrated in Figure~\ref{figrobustnesslocation}.
On the $x$-axis, the observation $\omega$ moves from 0 to 100. The
estimators $\hat{\mu}$ (left graph) and $\hat{\sigma}$ (right
graph) lie on the $y$-axis.

The influence of the outlier on a nonrobust inference is clearly
visible in the normal model (dashed lines) by the estimators
growing indefinitely as the outlier increases. For $\omega=100$, we
find $\hat{\mu}=4.55$ and $\hat{\sigma}=21.65$.
Using the normal quantile of 1.96,
this model thus suggests that 95\% of the observations should be
between $-$37.9 and 47.0, which is barely supported by data located
between $-$10 and 10,
and not at all by the outlier $\omega=100$.

Whole robustness is illustrated by the log-Pareto-tailed normal model
(solid lines).
We can see in Figure~\ref{figrobustnesslocation} that $\omega$
reaches its maximum influence around 16, where $\hat{\mu}$ and $\hat
{\sigma}$ are
approximately equal to 0.8 and 7. The influence of $\omega$ then
begins to decrease after this threshold
as $\hat{\mu}$ and $\hat{\sigma}$ eventually converge to their
corresponding MLE considering only
the nonoutlying observations $\mathbf{x_k}$, given by $\hat{\mu}=0$
and $\hat{\sigma}=6.06$.
For $\omega=100$, we find $\hat{\mu}=0.05$ and $\hat{\sigma}=6.28$.
Using the normal quantile of 1.96 (remember that this model is a
standard normal density except for the 2.5\%
log-Pareto tails),
this model thus suggests that 95\% of the observations should be
between $-$12.3 and 12.4, which is wholly supported by data, if
$\omega=100$ is considered as an outlier generated from the log-Pareto tails.

Finally, partial robustness is illustrated by the Student's $t$ model
(dotted-dashed lines).
For $\omega=100$, we find $\hat{\mu}=0.35$ and $\hat{\sigma}=8.44$.
Using the appropriate quantile of 2.147,
this model thus suggests that 95\% of the observations should be
between $-$17.8 and 18.5, which is partially supported by data located
between $-$10 and 10.
Note that as $\omega$ continues to grow beyond 100, our calculations
show that $\hat{\mu}$ decreases toward 0, and $\hat{\sigma}$
continues to grow
toward an upper limit of 8.71. This indicates that location estimation
using the Student's $t$ is wholly robust.
However, scale parameter estimation
is only partially robust, in the sense that the inference is
contaminated by the outlier, but only to a certain extent.

\subsection{Performance and simulations}\label{sec-performance}

We present here a brief study of the performance of the three models
described above (the robust log-Pareto-tailed normal,
the partially robust Student's $t$ and the popular but nonrobust
normal distributions) under three scenarios of simulations.
For each scenario and model, a sample of size $n=30$ is simulated
25,000 times, and the location and scale
parameters are estimated each time using the MLE.
Note that again, results using the Bayesian marginal posterior median
are very similar.
The performance is then measured by the mean square error (MSE). For
each scenario, the true values are $\mu=0$ and $\sigma=1$.
The MSE for the estimation of $\mu$ and $\sigma$ are
given in Tables~\ref{table1} and \ref{table2}, respectively.

In the first scenario, the samples are simulated from a $N(0,1)$.
We see that in the absence of outliers, the three models obtain the
same excellent performance both for the estimation
of the location (MSE${}=0.03$) and the scale (MSE${}=0.02$). This is rather
predictable, because the three
densities are very similar, if not identical, except for the tails. The
impact of the tails on the estimation is felt mainly in the presence of
outliers.

In the second scenario, we consider a mixture of normal distributions,
where an observation has a 90\% probability of being generated from a
$N(0,1)$ and 10\% from a $N(0,6)$.
A mixture of normal distributions is also studied in the third
scenario, where on average 95\% of the observations are generated from
a $N(0,1)$ and the remaining 5\%
from a $N(8,1)$.

%
\begin{table}[t]
\tabcolsep=0pt
\caption{Mean square error for MLE of $\mu$ under different scenarios ($n=30$)}\label{table1}
\begin{tabular*}{\tablewidth}{@{\extracolsep{\fill}}@{}lccc@{}}
\hline
&\multicolumn{3}{c@{}}{\textbf{Scenario}}\\[-6pt]
&\multicolumn{3}{c@{}}{\hrulefill}\\
\textbf{Model} & \textbf{100\% $\bolds{N(0,1)}$} & \textbf{10\% $\bolds{N(0,6)}$} & \textbf{5\% $\bolds{N(8,1)}$}\\
\hline
Log-Pareto-tailed normal & 0.03 & 0.05 & 0.07 \\
Student's $t$ & 0.03 & 0.06& 0.09 \\
Normal & 0.03& 0.15 & 0.29 \\
\hline
\end{tabular*}
\end{table}

%
\begin{table}[b]
\tabcolsep=0pt
\caption{Mean square error for MLE of $\sigma$ under different scenarios ($n=30$)}\label{table2}
\begin{tabular*}{\tablewidth}{@{\extracolsep{\fill}}@{}lccc@{}}
\hline
&\multicolumn{3}{c@{}}{\textbf{Scenario}}\\[-6pt]
&\multicolumn{3}{c@{}}{\hrulefill}\\
\textbf{Model} & \textbf{100\% $\bolds{N(0,1)}$} & \textbf{10\% $\bolds{N(0,6)}$} & \textbf{5\% $\bolds{N(8,1)}$}\\
\hline
Log-Pareto-tailed normal & 0.02 & 0.11& 0.09\\
Student's $t$ & 0.02 & 0.32& 0.30\\
Normal & 0.02& 1.46 & 1.14\\
\hline
\end{tabular*}
\end{table}

As for the estimation of $\mu$, we can see in Table~\ref{table1} that
both log-Pareto-tailed normal and
Student's $t$ models give very similar MSE for the two contaminated
scenarios (0.05 to 0.09), slightly larger than those of the 100\%
$N(0,1)$ scenario without outliers.
However, the normal model is clearly affected by the outliers as its
MSE increases to 0.15 and 0.29, respectively, for the second and third
scenarios.

The picture for the estimation of $\sigma$ is a bit different, as can
be seen in Table~\ref{table2}. For both scenarios,
the performance of the three models can be markedly discriminated in
accordance with known theory. The MSE are around 0.10 for the robust
log-Pareto normal model, around 0.30 for the partially robust Student
$t$ model and above 1
for the nonrobust normal model.

\section{Conclusion}\label{sec-conclusion}

Complete rejection of outliers has been investigated in a
location--scale parameter model.
The analysis has been done primarily in a Bayesian context, but it has
been extended to the frequentist approach with maximum likelihood estimators.
Essentially, asymptotic robustness is guaranteed if modeling is done
using a log-regularly varying distribution (with logarithmic tail
decay) and if $k>\max(l,r)$,
that is, if the number of nonoutliers
is larger than both the number of outliers at $-\infty$ and at
$+\infty$.
The first condition is easy to verify because it involves only the
tails of a density through a limit;
there are no integrals, derivatives or distribution functions involved.
The second condition is quite reasonable and intuitive.

We obtain the uniform convergence of the posterior density given the
complete sample to the density considering only the nonoutlying observations.
We also obtain the convergence in $L_1$, the convergence in
distribution, as well as the uniform convergence of the likelihoods.
Therefore any estimation of the location and scale parameters based on
posterior quantiles
or the maximum likelihood estimates is robust to outliers.

Even if the results are asymptotic, they are still useful in practice
with data, as illustrated by the threshold feature in Section~\ref
{sec-threshold}.
When one observation moves away from the rest of data, its influence on
the inference begins to increase gradually, because it brings
additional information
that helps us discriminate among the possible values of the parameter.
However, there comes a point where this moving observation conflicts
with the rest
of data. When this threshold is reached, the model automatically
resolves the conflict by progressively reducing the influence of the
outlying observation.
As the conflict grows infinitely, the impact of the outlier completely
disappears. This built-in feature is attractive in practice
in that conflict is managed in a sensitive and automatic way.

Estimating the location and scale parameters is common in statistics,
using, for instance, the well-known sample mean and standard deviation.
Results found in this paper can be readily used in practice to address
this problem
in a robust way, whether one prefers the Bayesian approach or maximum
likelihood estimation.
We consider a realistic sample of any size with multiple possible
outliers in any
direction. The assumption of a symmetric density $f$ with the same tail
behavior seems reasonable
for most of the applications. Because we do not know beforehand which
observations are going to be outlying,
it is generally desirable to give each density and each tail the same
weight, and to let the largest group
dominate in case of conflict.
The choice of the appropriate density is addressed in a practical way
by introducing the family of log-Pareto-tailed symmetric distributions.
Furthermore, the model allows us to add prior information on the
location and scale through a very general joint prior density, which
includes the possibility to choose a noninformative prior.

This paper can be generalized in different ways. For instance,
we can consider asymmetric densities $f$ with different tail behavior.
The family of log-regularly varying distributions could be widened to
consider, for instance,
distributions with a right tail proportional to $(1/z)\exp(-\delta
(\log z)^\gamma)$,
with $0<\gamma<1$ and $\delta>0$, which is an exponential
transformation of the function $\exp(-\delta z^\gamma)$.
Robustness to misspecification of the prior can also be investigated.

\section{Proofs}\label{sec-proof}

The proof of Proposition~\ref{proposition-proper} is given in
Section~\ref{proof-proposition-proper}, and the proof of result (a) of
Theorem~\ref{teo-main}
is given in Section~\ref{sec-proof-a}.

\subsection{Proof of Proposition~\texorpdfstring{\protect\ref
{proposition-proper}}{5}}\label{proof-proposition-proper}

To prove that $\pi(\mu,\sigma\mid\mathbf{x_n})$ is proper [the
proof for $\pi(\mu,\sigma\mid\mathbf{x_k})$ is omitted because it
is similar], it suffices to show that
the marginal $m(\mathbf{x_n})$ is finite.
Without loss of generality, we assume for convenience that
$x_1<x_2<\cdots<x_n$.
We also define the constant $\delta>0$ as half the minimum distance
between two observations, that is,
\[
\delta=\min_{i\in\{1,\ldots,n-1\}} \bigl\{(x_{i+1}-x_i)/2
\bigr\}.
\]

We first consider $\mu\in\re$ and $\delta/M\le\sigma<\infty$,
where $M$ is the constant of monotonicity given in equation~(\ref
{eqn-monotonic}).
Then we have
\begin{eqnarray*}
&& \int_{\delta/M}^{\infty}\int_{-\infty}^{\infty}
\pi(\mu,\sigma) \prod_{i=1}^{n} (1/\sigma)f
\bigl((x_i-\mu)/\sigma \bigr) \,d\mu\,d\sigma
\\
&&\qquad \za{\le}B^{n}\int_{\delta/M}^{\infty}(1/
\sigma)^{n}\int_{-\infty}^{\infty}(1/\sigma)f
\bigl((x_1-\mu)/\sigma \bigr) \,d\mu\,d\sigma
\\
&&\qquad \zb{=} B^{n}\int_{\delta/M}^{\infty}(1/
\sigma)^{n} \,d\sigma\int_{-\infty}^{\infty}f
\bigl( \mu' \bigr) \,d\mu'
\\
&&\qquad \zc{=} B^{n}(M/\delta)^{n-1}/(n-1)<\infty.
\end{eqnarray*}
In step $a$, we bound $\sigma\pi(\mu,\sigma)$ and $n-1$ densities
$f$ by $B$, where $B$ is given in (\ref{eqn-B}). In step $b$, we use
the change of variable $\mu'=(x_1-\mu)/\sigma$.
In step $c$, we use $n\ge2$ as assumed in the Bayesian context given
in Section~\ref{sec-model}.

We now consider $(x_{j-1}+x_{j})/2\le\mu\le(x_{j}+x_{j+1})/2$, for
$j=1,\ldots,n$ and $0<\sigma\le\delta/M$.
If we define $x_0:=-\infty$ and $x_{n+1}:=\infty$, the union of these
$n$ mutually disjoint intervals constitutes the real line, that is,
$-\infty<\mu<\infty$.
Then we have
\begin{eqnarray*}
&& \pi(\mu,\sigma)\prod_{i=1}^{n} (1/\sigma)f
\bigl((x_i-\mu)/\sigma \bigr)
\\
&&\qquad \za{\le}(1/\sigma)B\prod_{i=1}^{n}(1/
\sigma)f \bigl((x_i-\mu)/\sigma \bigr)
\\
&&\qquad =(1/\sigma)B f \bigl((x_j-\mu)/\sigma \bigr)\times(1/\sigma)\prod
_{i=1~(i\ne
j)}^{n}(1/\sigma)f \bigl((x_i-
\mu)/\sigma \bigr)
\\
&&\qquad \zb{\le}(1/\sigma)B f \bigl((x_j-\mu)/\sigma \bigr)\times (1/
\sigma) \bigl[(1/\sigma)f(\delta/\sigma) \bigr]^{n-1}
\\
&&\qquad \zc{\le}B(B/\delta)^{n-2}(1/\sigma) f \bigl((x_j-\mu)/
\sigma \bigr)\times(1/\sigma)^2 f(\delta/\sigma)
\\
&&\qquad \propto(1/\sigma) f \bigl((x_j-\mu)/\sigma \bigr)\times \bigl(
\delta/ \sigma^2 \bigr) f(\delta/\sigma).
\end{eqnarray*}
In step $a$, we bound $\sigma\pi(\mu,\sigma)$ by $B$.
In step $b$, we use $f((x_i-\mu)/\sigma)\le f(\delta/\sigma)$ by
the monotonicity of the tails of $f(z)$ since
$\llvert x_i-\mu\rrvert/\sigma\ge\delta/\sigma\ge\delta
(M/\delta)=M$,
because if $i\ne j$, we have
\[
\llvert x_i-\mu\rrvert\ge\min \bigl\{ (x_{j}-x_{j-1})/2,(x_{j+1}-x_{j})/2
\bigr\} \ge\delta.
\]
In step $c$, we bound $(1/\sigma)f(\delta/\sigma)$ by $B/\delta$
for $n-2$ terms.
Finally, we have
\begin{eqnarray*}
&& \int_{0}^{\delta/M} \bigl(\delta/\sigma^2
\bigr)f(\delta/\sigma)\int_{(x_{j-1}+x_{j})/2}^{(x_{j}+x_{j+1})/2} (1/\sigma)f
\bigl((x_j-\mu)/\sigma \bigr) \,d\mu\,d\sigma
\\
&& \qquad \le\int_{0}^{\infty}f \bigl(\sigma'
\bigr) \,d\sigma'\int_{-\infty}^{\infty
}f \bigl(
\mu' \bigr) \,d\mu'=1/2<\infty,
\end{eqnarray*}
where we use the changes of variable $\sigma'=\delta/\sigma$ and
$\mu'=(x_j-\mu)/\sigma$.

\subsection{Proof of result \textup{(a)} of Theorem~\texorpdfstring{\protect\ref
{teo-main}}{1}}\label{sec-proof-a}

Consider the model described in Section~\ref{sec-model},
and assume that the conditions of regularity on $f$ are satisfied. We
also assume that
$z f(z)\in L_{\rho}(\infty)$ and $k>\max(l,r)$, as given in
Theorem~\ref{teo-main}.
Two lemmas are first given, and the proof of result (a) follows.

%
\begin{lem}\label{cor-location--scale-transformation}
$\forall\lambda\ge0$, $\forall\tau\ge1$, there exists a constant
$D(\lambda,\tau)\ge1$ such that $z\in\re$ and
$(\mu,\sigma)\in[-\lambda,\lambda]\times[1/\tau, \tau
]\Rightarrow$
\[
1/D(\lambda,\tau)\le(1/\sigma) f \bigl((z-\mu)/\sigma \bigr)/f(z)\le D(\lambda,
\tau).
\]
\end{lem}

\begin{pf}
Proposition~\ref{prop-location--scale-transformation} states that
$(1/\sigma) f((z-\mu)/\sigma)/f(z)$ converges to 1 uniformly in any set
$(\mu,\sigma)\in E_{\lambda,\tau}$ as $z\rightarrow\infty$, where
$E_{\lambda,\tau}=[-\lambda,\lambda]\times[1/\tau, \tau]$.
Hence, $\forall\lambda\ge0$ and $\forall\tau\ge1$,
the ratio $(1/\sigma) f((z-\mu)/\sigma)/f(z)$ can be bounded, say by
$1/1.01$ and $1.01$,
if $\llvert z\rrvert$ is larger than a certain constant, say
$A(\lambda,\tau)$,
using the symmetry of $f$. Therefore, we choose $D(\lambda,\tau)\ge1.01$.

If $-A(\lambda,\tau)\le z\le A(\lambda,\tau)$, we observe that
$\llvert z-\mu\rrvert/\sigma$ is also bounded on $(\mu,\sigma
)\in E_{\lambda,\tau}$.
Therefore, since $f$ is continuous and strictly positive on $\re$, it
follows that $\forall\lambda\ge0$ and $\forall\tau\ge1$, we can
find a constant $D(\lambda,\tau)\ge1.01$
as large as we want such that the ratio $(1/\sigma) f((z-\mu)/\sigma
)/f(z)$ is bounded below by $1/D(\lambda,\tau)$ and above by
$D(\lambda,\tau)$, for any
$(\mu,\sigma)\in E_{\lambda,\tau}$.
\end{pf}

%
\begin{lem}\label{lem-convolution1}
There exists a constant $C>0$ such that
\[
\llvert z\rrvert\ge2 M\quad\Longrightarrow\quad\sup_{\mu\in\re}\frac{f(\mu
)f(z-\mu
)}{f(z)}
\le C,
\]
where $M$ is given in equation~(\ref{eqn-monotonic}).
\end{lem}
\begin{pf}
Let the constant $C=2 D(0,2) B$,
where $B$ is given in equation~(\ref{eqn-B}), and $D(0,2)$ comes from
Lemma~\ref{cor-location--scale-transformation}.
Consider $\llvert z\rrvert\ge2 M$.

First, consider $0\le\llvert\mu\rrvert\le\llvert z\rrvert/2$. We have
\[
\frac{f(\mu)f(z-\mu)}{f(z)}\za{\le} \frac{f(\mu)f(z/2)}{f(z)} \zb{\le} 2 D(0,2) f(\mu)\zc{\le}2
D(0,2) B=C. %
\]
In step $a$, we use $f(z-\mu)\le f(z/2)$ by the monotonicity of the
tails of $f$ since $\llvert z-\mu\rrvert\ge\llvert z\rrvert/2\ge
(2M)/2= M$.
In step $b$, we use $(1/2)f(z/2)/f(z)\le D(0,2)$ using Lemma~\ref
{cor-location--scale-transformation}.
In step $c$, we bound $f$ by $B$.

Second, consider $\llvert z\rrvert/2\le\llvert\mu\rrvert<\infty
$. We have
\[
\frac{f(\mu)f(z-\mu)}{f(z)} \le\frac{f(z/2)f(z-\mu)}{f(z)}\le2 D(0,2) f(z-\mu)\le2 D(0,2) B=C,
\]
using $f(\mu)\le f(z/2)$ in the first inequality by the monotonicity
of the tails of $f$ since $\llvert\mu\rrvert\ge\llvert
z\rrvert/2\ge(2M)/2= M$
and the same arguments as above for the other inequalities.
\end{pf}

We first observe that
\begin{eqnarray*}
&& \frac{m(\mathbf{x_n})}{m(\mathbf{x_k})\prod_{i=1}^{n}[f(x_i)]^{l_i+r_i}}
\\
&&\qquad= \frac{m(\mathbf{x_n})}{m(\mathbf{x_k})\prod_{i=1}^{n}[f(x_i)]^{l_i+r_i}} \int_{-\infty}^{\infty} \int
_{0}^{\infty}\pi(\mu,\sigma\mid \mathbf{x_n})
\,d\sigma\,d\mu
\\
&&\qquad= \int_{-\infty}^{\infty}\int_{0}^{\infty}
\frac
{\pi(\mu,\sigma)\prod_{i=1}^{n}
[(1/\sigma)f((x_i-\mu)/\sigma)
]^{k_i+l_i+r_i}}{m(\mathbf{x_k})\prod_{i=1}^{n}[f(x_i)]^{l_i+r_i}} \, d\sigma\,d\mu
\\
&&\qquad= \int_{-\infty}^{\infty}\int_{0}^{\infty}
\pi(\mu,\sigma\mid\mathbf{x_k}) \prod_{i=1}^{n}
\biggl[\frac
{(1/\sigma)f((x_i-\mu)/\sigma)}{f(x_i)} \biggr]^{l_i+r_i} \,d\sigma\,d\mu.
\end{eqnarray*}
Therefore, we show that the last integral converges to 1 as $\omega
\rightarrow\infty$ to prove result~(a).
If we use Lebesgue's dominated convergence theorem to pass the limit
$\omega\rightarrow\infty$ inside the integral, we have
\begin{eqnarray*}
&& \lim_{\omega\rightarrow\infty}\int_{-\infty}^{\infty} \int
_{0}^{\infty}\pi(\mu,\sigma\mid\mathbf{x_k})
\prod_{i=1}^{n} \biggl[\frac{(1/\sigma)f((x_i-\mu)/\sigma
)}{f(x_i)}
\biggr]^{l_i+r_i} \,d\sigma\,d\mu
\\
&&\qquad= \int_{-\infty}^{\infty}\int_{0}^{\infty}
\lim_{\omega
\rightarrow\infty} \pi(\mu,\sigma\mid\mathbf{x_k})\prod
_{i=1}^{n} \biggl[\frac
{(1/\sigma)f((x_i-\mu)/\sigma)}{f(x_i)}
\biggr]^{l_i+r_i} \,d\sigma\,d\mu
\\
&&\qquad= \int_{-\infty}^{\infty}\int_{0}^{\infty}
\pi(\mu,\sigma\mid\mathbf{x_k}) \,d\sigma\,d\mu= 1,
\end{eqnarray*}
using Proposition~\ref{prop-location--scale-transformation} in the
second equality and Proposition~\ref{proposition-proper} in the last one.
Note that pointwise convergence is sufficient, for any value of $\mu
\in\re$ and $\sigma>0$, once the limit is passed inside the integral.

However, in order to use Lebesgue's dominated convergence theorem, we
need to show that
$\pi(\mu,\sigma\mid\mathbf{x_k}) \prod_{i=1}^{n} [(1/\sigma
)f((x_i-\mu)/\sigma)/f(x_i) ]^{l_i+r_i}$
is bounded, for any value of $\omega\ge\xo$, by an integrable
function of $\mu$ and $\sigma$ that does not depend on~$\omega$.
The constant $\xo$ can be chosen as large as we want,
and some minimum values for $\xo$ will be given throughout the proof.

To achieve this, we divide the domain of integration into four
quadrants delineated by the axes $\mu=0$ and $\sigma=1$.
Note that the proofs are only given for the two quadrants in the region
of $\mu\ge0$ because the proofs for $\mu<0$ are similar.

We choose the constant $\xo$ larger than a certain threshold such that
the ranking of the set $\{\llvert x_i\rrvert\dvtx l_i+r_i=1\}$
remain unchanged for
all $\omega\ge\xo$. Given that each observation $x_i$ can be written
as $x_i=a_i+b_i \omega$, with $b_i=0$ if $k_i=1$, $b_i<0$ if $l_i=1$
and $b_i>0$ if $r_i=1$,
the ranking is therefore primarily determined by the values of
$\llvert b_i\rrvert$.
Then, without loss of generality, we assume for convenience that
\[
\min_{i\dvtx l_i+r_i=1} \bigl\{\llvert b_i\rrvert \bigr\}=1
\quad\mbox{and}\quad\omega=\min_{i\dvtx l_i+r_i=1} \bigl\{\llvert
x_i\rrvert \bigr\}.
\]

If $l_i+r_i=1$, we can use Lemma~\ref
{cor-location--scale-transformation}, with $x_i=a_i+b_i \omega
=b_i(\omega+a_i/b_i)$ and $\llvert b_i\rrvert\ge1$, to
establish that the ratio
$f(x_i)/f(\omega)$ is bounded, precisely by
\[
1/D \bigl(\llvert a_i/b_i\rrvert,\llvert
b_i\rrvert \bigr)\le\llvert b_i\rrvert
f(x_i)/f(\omega)\le D \bigl(\llvert a_i/b_i
\rrvert,\llvert b_i\rrvert \bigr). %
\]

\begin{qua}\label{q1}
Consider $0\le\mu<\infty$ and $1\le\sigma
<\infty$.
We have
\begin{eqnarray*}
&& \pi(\mu, \sigma\mid\mathbf{x_k})
\prod_{i=1}^{n}
\biggl[\frac{(1/\sigma)f((x_i-\mu)/\sigma)}{f(x_i)} \biggr]^{l_i+r_i}
\\
&&\qquad  \propto \frac{\pi(\mu,\sigma)}{\sigma^{n}}\prod
_{i=1}^{n}\frac
{f((x_i-\mu)/\sigma)}{ [f(x_i) ]^{l_i+r_i}}
\\
&&\qquad \za{\le}\frac{B}{\sigma^{n+1}}\prod_{i=1}^{n}
\frac
{D(\llvert a_i\rrvert,1)f((b_i \omega-\mu)/\sigma)}{
[f(x_i)
]^{l_i+r_i}}
\\
&&\qquad \zb{\le}\frac{1}{[f(\omega)]^{l+r}}\frac{B}{\sigma^{n+1}}\prod_{i=1}^{n}
D \bigl(\llvert a_i\rrvert,1 \bigr)f \bigl((b_i \omega-
\mu)/\sigma \bigr) \bigl[\llvert b_i\rrvert D \bigl(\llvert
a_i/b_i\rrvert,\llvert b_i\rrvert \bigr)
\bigr]^{l_i+r_i}
\\
&&\qquad \propto\frac{1}{[f(\omega)]^{l+r}}\frac{1}{\sigma^{n+1}}\prod_{i=1}^{n}
f \bigl((b_i \omega-\mu)/\sigma \bigr)
\\
&&\qquad \zc{=}\frac{1}{[f(\omega)]^{l+r}}\frac{[f(\mu/\sigma)]^k}{\sigma
^{n+1}} \prod_{i=1}^{n}
\bigl[f \bigl((b_i \omega-\mu)/\sigma \bigr) \bigr]^{l_i+r_i}
\\
&&\qquad = \frac{(1/\sigma)f(\mu/\sigma)}{\sigma^{k-1/2}} \biggl[\frac
{\omega/\sigma}{\omega f(\omega)} \biggr]^{l+r}
\frac{[f(\mu/\sigma
)]^{k-1}}{\sigma^{1/2}} \prod_{i=1}^{n} \bigl[f
\bigl((b_i \omega-\mu)/\sigma \bigr) \bigr]^{l_i+r_i}.
\end{eqnarray*}
In step $a$, we use $x_i=a_i+b_i \omega$ and
\[
f \bigl((x_i-\mu)/\sigma \bigr)=f \bigl((b_i \omega-
\mu)/ \sigma+a_i/\sigma \bigr)\le D \bigl(\llvert a_i
\rrvert,1 \bigr)f \bigl((b_i \omega-\mu)/\sigma \bigr)
\]
using Lemma~\ref{cor-location--scale-transformation}
since $\llvert a_i/\sigma\rrvert\le\llvert a_i\rrvert$.
We also bound $\sigma\pi(\mu,\sigma
)$ by $B$.
In step $b$, we use $1/f(x_i)\le\llvert b_i\rrvert D(\llvert
a_i/b_i\rrvert,\llvert b_i\rrvert)/f(\omega)$.
In step $c$, we set $b_i=0$ if $k_i=1$ and we use $f(-\mu/\sigma
)=f(\mu/\sigma)$ by symmetry of $f$.

It suffices to show that
\[
\biggl[\frac{\omega/\sigma}{\omega f(\omega)} \biggr]^{l+r}\frac
{[f(\mu/\sigma)]^{k-1}}{\sigma^{1/2}} \prod
_{i=1}^{n} \bigl[f \bigl((b_i
\omega-\mu)/\sigma \bigr) \bigr]^{l_i+r_i}<\infty,
\]
since $(1/\sigma)^{k-1/2}(1/\sigma)f(\mu/\sigma)$ is an integrable
function on Quadrant~\ref{q1},
\begin{eqnarray*}
\int_{1}^{\infty}(1/\sigma)^{k-1/2}\int
_{0}^{\infty}(1/\sigma)f(\mu/\sigma) \,d\mu\,d\sigma&\le&
\int_{1}^{\infty}(1/\sigma)^{k-1/2} \,d\sigma
= \frac{1}{k-3/2}\le2,
\end{eqnarray*}
since $k\ge2$.
To achieve this, we split the region of $\sigma$ into three parts between
$1<\omega^{1/2}<\omega/(2M)<\infty$, where $M$ is defined in
equation~(\ref{eqn-monotonic}).
Note that since $\omega\ge\xo$, this is well defined if $\xo> \max
(1,(2M)^2)$.

Consider $0\le\mu<\infty$ and $\omega/(2 M)\le\sigma<\infty$.
Then we have
\begin{eqnarray*}
&& \biggl[\frac{\omega/\sigma}{\omega f(\omega)} \biggr]^{l+r}\frac
{[f(\mu/\sigma)]^{k-1}}{\sigma^{1/2}} \prod_{i=1}^{n} \bigl[f \bigl((b_i
\omega-\mu)/\sigma \bigr) \bigr]^{l_i+r_i}
\\
&&\qquad  \za{\le} B^{n-1} \biggl[
\frac{\omega/\sigma}{\omega f(\omega
)} \biggr]^{l+r}\frac{1}{\sigma^{1/2}}
\zb{\le}B^{n-1}(2M)^{l+r+1/2}\frac{(1/\omega
)^{1/2}}{[\omega f(\omega)]^{l+r}}
\\
&&\qquad \zc{
\le}B^{n-1}(2M)^{l+r+1/2}\frac{(1/\omega)^{1/2}}{{(\log\omega
)^{-(\rho+1)(l+r)}}}
\\
&&\qquad\zd{\le}B^{n-1}(2M)^{l+r+1/2} \bigl[2(\rho+1) (l+r)/ \ee
\bigr]^{(\rho+1)(l+r)}<\infty.
\end{eqnarray*}
In step $a$, we use $f(\cdot)\le B$.
In step $b$, we use $\omega/\sigma\le2M$ and $(1/\sigma)\le
(2M)/\omega$.
In step $c$, we use $\omega f(\omega)>(\log\omega)^{-\rho-1}$ if
$\omega\ge\xo\ge A(1)$, where $A(1)$ comes from
Proposition~\ref{prop-dominance}.
In step $d$, it is purely algebraic to show that the maximum of $(\log
\omega)^{\beta}/\omega^{1/2}$ is $(2\beta/\ee)^{\beta}$
for $\omega>1$ and $\beta>0$, where $\beta=(\rho+1)(l+r)$ in our equation.

Now consider the two other parts combined (we will split them in the
next step), that is, $0\le\mu<\infty$ and $1\le\sigma\le\omega/(2M)$.
We have
\begin{eqnarray*}
\hspace*{-4pt}&& \biggl[\frac{\omega/\sigma}{\omega f(\omega)} \biggr]^{l+r}\frac
{[f(\mu/\sigma)]^{k-1}}{\sigma^{1/2}} \prod
_{i=1}^{n} \bigl[f \bigl((b_i
\omega-\mu)/\sigma \bigr) \bigr]^{l_i+r_i}
\\
\hspace*{-4pt}&&\qquad\za{\le} \biggl[\frac{\omega/\sigma}{\omega f(\omega
)} \biggr]^{l+r}
\frac{[f(\mu/\sigma)]^{k-1}}{\sigma^{1/2}} \prod_{i=1}^{n}
\bigl[f(b_i \omega/\sigma) \bigr]^{l_i} \bigl[f
\bigl((b_i \omega-\mu)/\sigma \bigr) \bigr]^{r_i}
\\
\hspace*{-4pt}&&\qquad=\frac{[f(\mu/\sigma)]^{k-r-1}}{\sigma^{1/2}}
\prod_{i=1}^{n}
\biggl[\frac{(\omega/\sigma)f(b_i \omega/\sigma
)}{\omega f(\omega)} \biggr]^{l_i+r_i} \biggl[\frac{f(b_i \omega
/\sigma-\mu/\sigma)f(\mu/\sigma)}{f(b_i
\omega/\sigma)}
\biggr]^{r_i}
\\
\hspace*{-4pt}&&\qquad\zb{\le}B^{k-r-1}C^r \frac{1}{\sigma^{1/2}} \prod
_{i=1}^{n} \biggl[\frac{(\omega/\sigma)f(b_i \omega/\sigma
)}{\omega f(\omega)}
\biggr]^{l_i+r_i}
\\
\hspace*{-4pt}&&\qquad\zc{\le}B^{k-r-1}C^r \frac{1}{\sigma^{1/2}} \biggl[
\frac{(\omega/\sigma)f(\omega/\sigma)}{\omega f(\omega
)} \biggr]^{l+r}.
\end{eqnarray*}
In step $a$, we use $f((b_i \omega-\mu)/\sigma)\le f(b_i \omega
/\sigma)$ if $l_i=1$ (which means $b_i<0$) by the monotonicity of the
tails of $f$ since
$\llvert b_i \omega-\mu\rrvert/\sigma\ge\llvert b_i\rrvert
\omega/\sigma\ge\llvert b_i\rrvert (2M) \ge
2M \ge M$. In step $b$, we use $f(\mu/\sigma)\le B$ and
we use Lemma~\ref{lem-convolution1} since $\llvert b_i\rrvert
\omega/\sigma\ge
\llvert b_i\rrvert (2M) \ge2M$. In step $c$, we use
$f(b_i \omega/\sigma)\le f(\omega/\sigma)$ by the monotonicity of
the tails of $f$ since $\llvert b_i\rrvert\omega/\sigma\ge\omega
/\sigma\ge2M
\ge M$.\vspace*{1pt}

Consider $0\le\mu<\infty$ and $\omega^{1/2}\le\sigma\le\omega/(2M)$.
We have
\begin{eqnarray*}
\frac{1}{\sigma^{1/2}} \biggl[\frac{(\omega/\sigma)f(\omega
/\sigma)}{\omega f(\omega
)} \biggr]^{l+r}& \za{\le}&
B^{l+r} \frac{(1/\omega)^{1/4}}{[\omega f(\omega)]^{l+r}}
\\
&\zb {\le}& B^{l+r} \frac{(1/\omega)^{1/4}}{(\log\omega)^{-(\rho+1)(l+r)}}
\\
&\zc{\le} & B^{l+r} \bigl[4(\rho+1) (l+r)/\ee \bigr]^{(\rho
+1)(l+r)}<
\infty.
\end{eqnarray*}
In step $a$, we use $(\omega/\sigma)f(\omega/\sigma)\le B$ and
$(1/\sigma)^{1/2}\le(1/\omega)^{1/4}$.
In step $b$, we use $\omega f(\omega)>(\log\omega)^{-\rho-1}$ if
$\omega\ge\xo\ge A(1)$, where $A(1)$ comes from
Proposition~\ref{prop-dominance}.
In step $c$, it is purely algebraic to show that the maximum of $(\log
\omega)^{\beta}/\omega^{1/4}$ is $(4\beta/\ee)^{\beta}$
for $\omega>1$ and $\beta>0$, where $\beta=(\rho+1)(l+r)$ in our equation.

Finally consider $0\le\mu<\infty$ and $1\le\sigma\le\omega^{1/2}$.
Then we have
\[
\frac{1}{\sigma^{1/2}} \biggl[\frac{(\omega/\sigma)f(\omega
/\sigma)}{\omega f(\omega
)} \biggr]^{l+r}\za{\le}
\biggl[\frac{\omega^{1/2}f(\omega^{1/2})}{\omega f(\omega)} \biggr]^{l+r}\zb{\le} 2^{(\rho+1)(l+r)}<\infty.
\]
In step $a$, we use $1/\sigma\le1$, and we use $(\omega/\sigma
)f(\omega/\sigma)\le\omega^{1/2} f(\omega^{1/2})$ by the
monotonicity of the tails of
$\llvert z\rrvert f(z)$ since $\omega/\sigma\ge\omega^{1/2}\ge
\xo^{1/2}\ge M$
if $\xo\ge M^2$. In step~$b$, we use $\omega^{1/2} f(\omega
^{1/2})/(\omega
f(\omega))\le2(1/2)^{-\rho}=2^{\rho+1}$
if $\omega\ge\xo\ge A(1,2)$, where $A(1,2)$ comes from the
definition of a log-regularly varying function.
\end{qua}

\begin{qua}\label{q2}
Consider $-\infty<\mu\le0$ and $1\le\sigma
<\infty$.
The proof for Quadrant~\ref{q2} is similar to that of Quadrant~\ref{q1}.
\end{qua}

\begin{qua}\label{q3}
Consider $-\infty<\mu\le0$ and $0<\sigma\le1$.
The proof for Quadrant~\ref{q3} is similar to that of Quadrant~\ref{q4}, given below.
The condition
$k>r$ is therefore replaced by $k>l$. Note that $k>\max(l,r)$ is
assumed in Theorem~\ref{teo-main}.
\end{qua}

\begin{qua}\label{q4}
Consider $0\le\mu<\infty$ and $0<\sigma\le1$.
We need to show, actually, that
\begin{eqnarray*}
&& \lim_{\omega\rightarrow\infty} \int_{0}^{\infty}\int
_{0}^{1} \pi(\mu,\sigma\mid\mathbf{x_k})
\prod_{i=1}^{n} \biggl[\frac
{(1/\sigma)f((x_i-\mu)/\sigma)}{f(x_i)}
\biggr]^{l_i+r_i} \,d\sigma\,d\mu
\\
&&\qquad = \int_{0}^{\infty}\int_{0}^{1}
\pi(\mu,\sigma\mid\mathbf{x_k}) \,d\sigma\,d\mu.
\end{eqnarray*}
For Quadrant~\ref{q1}, we show this result when we integrate $\sigma$ between
1 and $\infty$. We bound the integrand of the left term, for any value
of $\omega\ge\xo$, by an integrable
function of $\mu$ and $\sigma$ that does not depend on $\omega$,
in order to use Lebesgue's dominated convergence theorem to pass the
limit $\omega\rightarrow\infty$ inside the integral.
For Quadrant~\ref{q4}, we proceed slightly differently. We begin by breaking
down the left term into two parts as follows:
\begin{eqnarray*}
&&\lim_{\omega\rightarrow\infty} \int_{0}^{\infty}\int
_{0}^{1} \pi(\mu,\sigma\mid\mathbf{x_k})
\prod_{i=1}^{n} \biggl[\frac
{(1/\sigma)f((x_i-\mu)/\sigma)}{f(x_i)}
\biggr]^{l_i+r_i} \,d\sigma\,d\mu
\\
&&\qquad=\lim_{\omega\rightarrow\infty} \int_{0}^{\infty}
\int_{0}^{1} \pi(\mu,\sigma\mid
\mathbf{x_k}) \prod_{i=1}^{n}
\biggl[\frac
{(1/\sigma)f((x_i-\mu)/\sigma)}{f(x_i)} \biggr]^{l_i+r_i}
\\
&&\hspace*{161pt}{}\times \mathbh{1}_{[0,\omega/2]}(\mu)
\,d\sigma\,d\mu
\\
&&\quad\qquad{} +\lim_{\omega\rightarrow\infty} \int_{\omega
/2}^{\infty}
\int_{0}^{1} \pi(\mu,\sigma\mid
\mathbf{x_k}) \prod_{i=1}^{n}
\biggl[\frac
{(1/\sigma)f((x_i-\mu)/\sigma)}{f(x_i)} \biggr]^{l_i+r_i} \,d\sigma\,d\mu,
\end{eqnarray*}
where the indicator function $\mathbh{1}_{A}(\mu)$ is equal to 1 if
$\mu\in A$, and equal to 0 otherwise. We then show that
the first part is equal to the integral $\int_{0}^{\infty}\int_{0}^{1}\pi(\mu,\sigma\mid\mathbf{x_k}) \,d\sigma\,d\mu$,
and the second part is equal to 0.

For the first equality, we again use Lebesgue's dominated convergence
theorem to pass the limit $\omega\rightarrow\infty$ inside the
integral. We have
\begin{eqnarray*}
&&\lim_{\omega\rightarrow\infty}\int_{0}^{\infty}\int
_{0}^{1}\pi(\mu,\sigma\mid\mathbf{x_k})
\prod_{i=1}^{n} \biggl[\frac
{(1/\sigma)f((x_i-\mu)/\sigma)}{f(x_i)}
\biggr]^{l_i+r_i} \mathbh{1}_{[0,\omega/2]}(\mu) \,d\sigma\,d\mu
\\
&&\qquad= \int_{0}^{\infty}\int_{0}^{1}
\pi(\mu,\sigma\mid\mathbf{x_k})\lim_{\omega\rightarrow\infty
}\prod
_{i=1}^{n} \biggl[\frac{(1/\sigma)f((x_i-\mu)/\sigma
)}{f(x_i)}
\biggr]^{l_i+r_i}
\\
&&\hspace*{162pt}{}\times \mathbh{1}_{[0,\omega/2]}(\mu) \,d\sigma\,d\mu
\\
&&\qquad =\int_{0}^{\infty}\int_{0}^{1}
\pi(\mu,\sigma\mid\mathbf{x_k})\times1 \times\mathbh
{1}_{[0,\infty)}(\mu) \,d\sigma\,d\mu
\\
&&\qquad =\int_{0}^{\infty}
\int_{0}^{1} \pi(\mu,\sigma\mid
\mathbf{x_k}) \,d\sigma\,d\mu,
\end{eqnarray*}
using Proposition~\ref{prop-location--scale-transformation} in the
second equality.
Note that pointwise convergence is sufficient, for any value of $\mu
\in\re$ and $\sigma>0$,
once the limit is passed inside the integral.
However, in order to use Lebesgue's dominated convergence theorem, we
need to show that
for any value of $\omega\ge\xo$, the integrand
is bounded by an integrable function of $\mu$ and $\sigma$ that does
not depend on $\omega$.

Consider $0\le\mu\le\omega/2$ (the integrand is equal to 0 if
$\omega/2< \mu<\infty$) and $0<\sigma\le1$. We have
\begin{eqnarray*}
&& \pi(\mu,\sigma\mid\mathbf{x_k}) \prod
_{i=1}^{n} \biggl[\frac{(1/\sigma)f((x_i-\mu)/\sigma
)}{f(x_i)}
\biggr]^{l_i+r_i}\mathbh{1}_{[0,\omega/2]}(\mu)
\\
&&\qquad\za{\le}\pi(\mu,\sigma\mid\mathbf{x_k}) \prod
_{i=1}^{n} \biggl[2D(0,2)\frac{(1/\sigma)f((x_i-\mu)/\sigma
)}{f(x_i/2)}
\biggr]^{l_i+r_i}
\\
&&\qquad\propto\pi(\mu,\sigma\mid\mathbf{x_k}) \prod
_{i=1}^{n} \biggl[\frac{(1/\sigma)f((x_i-\mu)/\sigma
)}{f(x_i/2)}
\biggr]^{l_i+r_i}
\\
&&\qquad\zb{\le}\pi(\mu,\sigma\mid\mathbf{x_k}) \prod
_{i=1}^{n} \biggl[\frac{f(x_i-\mu)}{f(x_i/2)}
\biggr]^{l_i+r_i}
\\
&&\qquad \zc{\le}\pi(\mu,\sigma\mid\mathbf{x_k}),
\end{eqnarray*}
and $\pi(\mu,\sigma\mid\mathbf{x_k})$ is an integrable function.
In step $a$, we use $\mathbh{1}_{[0,\omega/2]}(\mu)= 1$ and
$(1/2)f(x_i/2)/f(x_i)\le D(0,2)$ using Lemma~\ref
{cor-location--scale-transformation}. In step $b$, we use $(\llvert
x_i-\mu
\rrvert/\sigma)f((x_i-\mu)/\sigma)\le\llvert x_i-\mu\rrvert
f(x_i-\mu)$
by the monotonicity of the tails of
$\llvert z\rrvert f(z)$, and in step $c$ we use $f(x_i-\mu)\le
f(x_i/2)$ by the
monotonicity of the tails of
$f(z)$
since $\llvert x_i-\mu\rrvert/\sigma\ge\llvert x_i-\mu\rrvert
\ge\llvert x_i\rrvert/2\ge\omega/2 \ge\xo
/2\ge M$, if we choose $\xo\ge2 M$.
Note that the condition $\mu\le\omega/2 (\le x_i/2)$ is used only
to justify $\llvert x_i-\mu\rrvert\ge\llvert x_i\rrvert
/2$ when $r_i=1$.

Now we show the second equality, that is,
\[
\lim_{\omega\rightarrow\infty} \int_{\omega/2}^{\infty}\int
_{0}^{1} \pi(\mu,\sigma\mid\mathbf{x_k})
\prod_{i=1}^{n} \biggl[\frac
{(1/\sigma)f((x_i-\mu)/\sigma)}{f(x_i)}
\biggr]^{l_i+r_i} \,d\sigma\,d\mu=0.
\]
We first bound above the integrand, and then we show that the integral
of the upper bound converges to 0 as $\omega\rightarrow\infty$.

Consider $\omega/2\le\mu<\infty$ and $0<\sigma\le1$. We have
\begin{eqnarray*}
&& \pi(\mu,\sigma\mid\mathbf{x_k})\prod
_{i=1}^{n} \biggl[\frac
{(1/\sigma)f((x_i-\mu)/\sigma)}{f(x_i)}
\biggr]^{l_i+r_i}
\\
&&\qquad\za{\le} \bigl[2D(0,2) \bigr]^l\pi(\mu,\sigma\mid
\mathbf{x_k}) \prod_{i=1}^{n}
\biggl[\frac{\llvert b_i\rrvert D(\llvert a_i/b_i\rrvert,\llvert
b_i\rrvert)(1/\sigma)f((x_i-\mu)/\sigma
)}{f(\omega)} \biggr]^{r_i}
\\
&&\qquad\propto\pi(\mu,\sigma)\prod_{i=1}^{n}
\bigl[(1/\sigma)f \bigl((a_i-\mu)/\sigma \bigr)
\bigr]^{k_i} \biggl[\frac{(1/\sigma)f((x_i-\mu)/\sigma
)}{f(\omega)} \biggr]^{r_i}
\\
&&\qquad\zb{\le}(1/\sigma)B \bigl[4 D(0,4) (1/\sigma)f(\omega /\sigma)
\bigr]^{k}\prod_{i=1}^{n} \biggl[
\frac{(1/\sigma)f((x_i-\mu
)/\sigma)}{f(\omega)} \biggr]^{r_i}
\\
&&\qquad\propto(1/\sigma) \bigl[(1/\sigma)f(\omega/\sigma) \bigr]^{k}
\prod_{i=1}^{n} \biggl[\frac{(1/\sigma)f((x_i-\mu)/\sigma
)}{f(\omega)}
\biggr]^{r_i}
\\
&&\qquad\zc{\le} (1/\sigma) \bigl[(1/\sigma)f(\omega/\sigma )
\bigr]^{k-r}\prod_{i=1}^{n}
\bigl[(1/\sigma)f \bigl((x_i-\mu)/\sigma \bigr) \bigr]^{r_i}
\\
&&\qquad\zd{=} (1/\sigma) \bigl[(1/\sigma)f(\omega/\sigma) \bigr]^{k-r}
\prod_{i=1}^{r}(1/\sigma)f
\bigl((x_i- \mu)/\sigma \bigr).
\end{eqnarray*}
In step $a$, we use $(1/\sigma)f((x_i-\mu)/\sigma)/f(x_i)\le
2D(0,2)$ if $l_i=1$, using the same arguments given above for the case
$0\le\mu\le\omega/2$.
We also use $1/f(x_i)\le\llvert b_i\rrvert D(\llvert
a_i/b_i\rrvert,\llvert b_i\rrvert)/f(\omega)$ if $r_i=1$.
In step $b$, we bound $\sigma\pi(\mu,\sigma)$ by $B$. We also use
\[
f \bigl((a_i-\mu)/\sigma \bigr)\le f \bigl((1/4)\omega/\sigma
\bigr) \le4 D(0,4) f(\omega/\sigma)
\]
if $k_i=1$
using the monotonicity of the tails of $f(z)$ in the first inequality
since, if we define $a_{(k)}=\max_{i\dvtx k_i=1}\{\llvert a_i\rrvert
\}$ with $\omega
\ge\xo\ge4 a_{(k)}$, we have
$\llvert a_i-\mu\rrvert/\sigma=(\mu-a_i)/\sigma\ge(\omega
/2-a_{(k)})/\sigma
\ge(\omega/2-\omega/4)/\sigma=(1/4)\omega/\sigma
\ge\omega/4 \ge\xo/4 \ge M$ if we choose $\xo\ge4 M$. We use
Lemma~\ref{cor-location--scale-transformation} in the second inequality.
In step $c$, we use $(\omega/\sigma)f(\omega/\sigma)\le\omega
f(\omega)$, using the monotonicity of the tails of $\llvert z\rrvert f(z)$
since $\omega/\sigma\ge\omega\ge\xo\ge M$ if we choose $\xo\ge M$.
In step $d$, we assume for convenience and without loss of generality
that the right outliers are denoted by \mbox{$x_1<x_2<\cdots<x_r$}.

We now split the real line (which includes the region $\omega/2\le\mu
<\infty$) into $r$ mutually disjoint intervals given by
$(x_{j-1}+x_{j})/2\le\mu\le(x_{j}+x_{j+1})/2$, for $j=1,\ldots,r$,
where we define $x_0:=-\infty$ and $x_{r+1}:=\infty$.
We also define the constant $\delta>0$ as
\[
\delta=\min_{i\in\{1,\ldots,r-1\}} \bigl\{(x_{i+1}-x_i)/2
\bigr\}.
\]
Consider $(x_{j-1}+x_{j})/2\le\mu\le(x_{j}+x_{j+1})/2$, for
$j=1,\ldots,r$ and $0<\sigma\le1$.
Then we have
\begin{eqnarray*}
&& (1/\sigma) \bigl[(1/\sigma)f(\omega/\sigma) \bigr]^{k-r}\prod
_{i=1}^{r}(1/\sigma)f \bigl((x_i-
\mu)/ \sigma \bigr)
\\
&&\qquad\za{\le}(B/\delta)^{r-1} (1/\sigma) \bigl[(1/\sigma)f( \omega/
\sigma) \bigr]^{k-r}(1/\sigma)f \bigl((x_j-\mu)/\sigma
\bigr)
\\
&&\qquad\zb{\le} (B/\delta)^{r-1}B^{k-r-1}\omega ^{-(k-r)}
\bigl(\omega/\sigma^2 \bigr)f(\omega/\sigma)\times(1/ \sigma)f
\bigl((x_j-\mu)/\sigma \bigr).
\end{eqnarray*}
In step $a$, we use, for $i\ne j$, $(1/\sigma)f((x_i-\mu)/\sigma)\le
B/\llvert x_i-\mu\rrvert\le B/\delta$, where we bound $\llvert
z\rrvert f(z)$ by $B$, and we
use $\llvert x_i-\mu\rrvert\ge\delta$ because if $i\ne j$, we have
\[
\llvert x_i-\mu\rrvert\ge\min \bigl\{ (x_{j}-x_{j-1})/2,(x_{j+1}-x_{j})/2
\bigr\} \ge\delta.
\]
In step $b$, we use $(\omega/\sigma)f(\omega/\sigma)\le B$ for
$k-r-1$ terms.
Finally, we have
\begin{eqnarray*}
&& \omega^{-(k-r)}\int_{0}^{1} \bigl(
\omega/ \sigma^2 \bigr)f(\omega/\sigma)\int_{(x_{j-1}+x_{j})/2}^{(x_{j}+x_{j+1})/2}(1/
\sigma)f \bigl((x_j-\mu)/\sigma \bigr) \,d\mu\,d\sigma
\\
&&\qquad\le\omega^{-(k-r)}\int_{0}^{\infty}
\bigl( \omega/\sigma^2 \bigr)f(\omega/\sigma)\int_{-\infty}^{\infty}(1/
\sigma)f \bigl((x_j-\mu)/\sigma \bigr) \,d\mu\,d\sigma
\\
&&\qquad\za{=} \omega^{-(k-r)}\int_{0}^{\infty}f
\bigl(\sigma' \bigr) \,d\sigma'\int
_{-\infty}^{\infty}f \bigl(\mu' \bigr) \,d
\mu' \le\omega^{-(k-r)}\zb{ \rightarrow} 0\qquad\mbox{as }
\omega\rightarrow\infty.
\end{eqnarray*}
In step $a$, we use the changes of variable $\sigma'=\omega/\sigma$
and $\mu'=(x_j-\mu)/\sigma$. In step~$b$, we use the condition $k>r$.
\end{qua}

\section*{Acknowledgments}\label{sec-Acknowledgment}
We thank the Associate Editor and the referees for very helpful comments.




%

\printaddresses
\end{document}